\documentstyle[12pt]{article}
\begin{document}
\newtheorem{theorem}      {Th\'eor\`eme}[section]
\newtheorem{theorem*}     {theorem}
\newtheorem{proposition}  [theorem]{Proposition}
\newtheorem{definition}   [theorem]{Definition}
\newtheorem{e-lemme}        [theorem]{Lemma}
\newtheorem{cor}   [theorem]{Corollaire}
\newtheorem{resultat}     [theorem]{R\'esultat}
\newtheorem{eexercice}    [theorem]{Exercice}
\newtheorem{rrem}    [theorem]{Remarque}
\newtheorem{pprobleme}    [theorem]{Probl\`eme}
\newtheorem{eexemple}     [theorem]{Exemple}
\newcommand{\preuve}      {\paragraph{Preuve}}
\newenvironment{probleme} {\begin{pprobleme}\rm}{\end{pprobleme}}
\newenvironment{remarque} {\begin{rremarque}\rm}{\end{rremarque}}
\newenvironment{exercice} {\begin{eexercice}\rm}{\end{eexercice}}
\newenvironment{exemple}  {\begin{eexemple}\rm}{\end{eexemple}}
%
%
\newtheorem{e-theo}      [theorem]{Theorem}
\newtheorem{theo*}     [theorem]{Theorem}
\newtheorem{e-pro}  [theorem]{Proposition}
\newtheorem{e-def}   [theorem]{Definition}
\newtheorem{e-lem}        [theorem]{Lemma}
\newtheorem{e-cor}   [theorem]{Corollary}
\newtheorem{e-resultat}     [theorem]{Result}
\newtheorem{ex}    [theorem]{Exercise}
\newtheorem{e-rem}    [theorem]{Remark}
\newtheorem{prob}    [theorem]{Problem}
\newtheorem{example}     [theorem]{Example}
\newcommand{\proof}         {\paragraph{Proof~: }}
\newcommand{\hint}          {\paragraph{Hint}}
\newcommand{\heuristicproof}{\paragraph{heuristic proof}}
\newenvironment{e-probleme} {\begin{e-pprobleme}\rm}{\end{e-pprobleme}}
\newenvironment{e-remarque} {\begin{e-rremarque}\rm}{\end{e-rremarque}}
\newenvironment{e-exercice} {\begin{e-eexercice}\rm}{\end{e-eexercice}}
\newenvironment{e-exemple}  {\begin{e-eexemple}\rm}{\end{e-eexemple}}
\newcommand{\1}        {{\bf 1}}
\newcommand{\pp}       {{{\rm I\!\!\! P}}}
\newcommand{\qq}       {{{\rm I\!\!\! Q}}}
\newcommand{\B}        {{{\rm I\! B}}}
\newcommand{\cc}       {{{\rm I\!\!\! C}}}
\newcommand{\N}        {{{\rm I\! N}}}
\newcommand{\R}        {{{\rm I\! R}}}
\newcommand{\D}        {{{\rm I\! D}}}
\newcommand{\Z}       {{{\rm Z\!\! Z}}}
\newcommand{\C}        {{\bf C}}        
\newcommand{\rank}{\hbox{rank}}
\newcommand{\CC}{{\cal C}}
\def \Re {{\rm Re\,}}
\def \Im {{\rm Im\,}}
\def\Hom{{\rm Hom\,}}
\def\Lip{{\rm Lip}}
%
%
\newcommand{\dontforget}[1]
{{\mbox{}\\\noindent\rule{1cm}{2mm}\hfill don't forget : #1
\hfill\rule{1cm}{2mm}}\typeout{---------- don't forget : #1 ------------}}
\newcommand{\note}[2]
{ \noindent{\sf #1 \hfill \today}

\noindent\mbox{}\hrulefill\mbox{}
\begin{quote}\begin{quote}\sf #2\end{quote}\end{quote}
\noindent\mbox{}\hrulefill\mbox{}
\vspace{1cm}
}
\title{\huge Filling hypersurfaces by discs in  almost complex
manifolds of  dimension 2}
\author{ Alexandre Sukhov{*} and Alexander Tumanov{**}}
\date{}
\maketitle

{\small {*}Universit\'e des Sciences et Technologies de Lille,
Laboratoire Paul Painlev\'e,
U.F.R. de Math\'e-matique, 59655 Villeneuve d'Ascq, Cedex, France,
sukhov@math.univ-lille1.fr

{**} University of Illinois, Department of Mathematics
1409 West Green Street, Urbana, IL 61801, USA, tumanov@math.uiuc.edu}
\bigskip

Abstract. We study pseudoholomorphic discs with boundaries attached
to a real hypersurface $E$ in an almost complex manifold of
dimension 2. We prove that if $E$ contains no  discs, then they
fill a one sided neighborhood of $E$.
\bigskip


\bigskip

\section{Introduction}
The modern development of the analysis on almost complex manifolds
began after the work of Gromov \cite{Gr} who discovered its remarkable
applications  in  symplectic geometry. One of the main tools of his
approach involves the Bishop discs, which are pseudoholomorphic
discs with boundaries attached to a prescribed real
submanifold; they have been used later by many authors.

On the other hand, the method of Bishop discs is well known and
powerful in complex analysis of several variables
(see, e.g., \cite{Tr2,Tu2}).
The main goal of this paper is to develop a systematic
general approach to the theory of Bishop discs in almost
complex manifolds.
As an object of the study we choose a well known
nontrivial result for the standard complex structure
and generalize it to the almost complex category.
Our main result is the following.
\begin{e-theo}
\label{MainTheorem}
Let $E$ be a real hypersurface in an almost complex manifold $(M,J)$
of complex dimension $2$. Suppose that $E$ contains no $J$-holomorphic
discs. Then Bishop discs of $E$ fill a one-sided neighborhood of
every point of $E$.
\end{e-theo}

In the case of the standard complex structure, the corresponding
result was obtained by Tr\'epreau \cite{Tr}
for real hypersurfaces in $\cc^n$ and by the second named
author \cite{Tu} for generic submanifolds of
arbitrary codimension in $\cc^n$. Although the restriction on the
dimension is used in our paper, this is mainly due to  technical
reasons and  we believe that after a suitable modification our
approach can be generalized to higher dimension.

In fact, a stronger result is proved in \cite{Tr, Tu}:
if Bishop discs don't fill a one-sided neighborhood of $p\in E$,
then there exists a disc in E passing through $p$.
Thus we leave open an interesting question
whether the same can be done in the almost complex case.
We may want to address it elsewhere.

We now  explain  the organization of the paper and the main steps
of our approach. Section 2 contains basic properties of almost
complex manifolds, pseudoholomorphic discs and some special coordinate
systems on almost complex manifolds used throughout the paper.
In Section 3 we derive different versions of the Bishop equation
for pseudoholomorphic discs. It can be viewed as a non-linear
boundary value Riemann-Hilbert type problem of a quasilinear elliptic
system in the unit disc. We prove that the Bishop discs attached
to a real hypersurface at a given point form a Banach manifold and
give an efficient local parametrization of this manifold. As an
application in Section 4 we prove Theorem \ref{MainTheorem} under
some additional restrictions (pseudoconvexity or finite type
conditions). In particular, we prove Theorem \ref{MainTheorem}
in the real analytic category. Our considerations there are
primarily geometrical. They are based on non-isotropic dilations
of suitable local coordinates which allow to represent the above
mentioned Riemann-Hilbert type problem as a small perturbation
of the corresponding problem for the standard complex
structure in $\cc^n$.

Section 5 is dedicated to the analysis of the Bishop equation for
pseudoholomorphic discs and is principal for our approach.
We parametrize the tangent space to the manifold of Bishop discs
by a space of holomorphic functions and solve the corresponding
linearized boundary value problem by means of the generalized
Schwarz integral. It turns out that the proof of
Theorem \ref{MainTheorem} given in \cite{Tu} for the standard
complex structure does not go through in the almost complex category.
The proof in \cite{Tu} is based on the notion of the defect of a disc.
For a Bishop disc through a fixed point of a real manifold $E$,
the infinitesimal perturbations of the direction of the disc at
the fixed point and those of another boundary point of the disc
are restricted in the same way, according to the defect of the disc.
We were unable to find an analogue of this phenomenon in our
case of complex dimension 2. In higher dimension no such analogue
is possible due to an example by Ivashkovich and Rosay \cite{ivro},
in which all Bishop discs through a fixed point of a hypersurface
$E$ lie in $E$ and cover all of $E$.
Baouendi, Rothschild, and Tr\'epreau \cite{BaTrRo} interpret the
defect of a disc in terms of its lifts attached to the conormal
bundle of $E$ in the cotangent bundle $T^*\cc^n$.
Although an almost complex structure admits natural lifts to the
cotangent bundle of an almost complex manifold (see, e.g., \cite{YI}),
seemingly, they don't give rise to a correct notion of the defect
of Bishop discs.

In this paper we develop a new approach, in particular we give a
different proof of the main result of \cite{Tu} for a hypersurface
in $\cc^2$. Our key result is the following.

\begin{e-theo}
\label{keytheo}
Let  $E$ be a real hypersurface in an almost complex manifold $(M,J)$ of complex
dimension $2$ and let $f_0$ be a small enough embedded Bishop disc attached to $E$
at a point $p \in E$ (that is $f_0(1) = p$) and tangent to $E$ at $p$.
Suppose that every Bishop disc $f$ attached at $p$ and close enough to $f_0$
also is tangent to $E$ at $p$. Then the Levi form of $E$ vanishes identically
along the boundary of $f_0$.
\end{e-theo}

For the standard complex structure, this result follows from \cite{Tu}.
In the almost complex case the result is new.
We show that the condition of tangency to the hypersurface $E$ of all Bishop discs attached to $E$
at a point $p$ is equivalent to the vanishing of some non-linear
operator $\Phi$ defined on
the Banach space of discs and valued in the space of smooth functions on the unit circle.
We show that the Frechet derivative $\dot \Phi$ of $\Phi$ up to smoother terms is equal to the multiplication
operator by the Levi determinant of $E$ and this allows to conclude the proof
of Theorem \ref{keytheo}. Now if in the hypothesis of Theorem
\ref{MainTheorem} $E$ does not admit a transversal Bishop disc attached at a  given
point $p$ then Theorem \ref{keytheo} implies that $E$ is Levi flat along the boundaries of the Bishop
discs through $p$ which allows to construct holomorphic discs
contained in $E$ and  leads to a contradiction. Hence there exists a transversal Bishop disc,
whose perturbations fill a one sided neighborhood of $E$.
Theorem \ref{keytheo} (and the method of its proof) gives a new
powerful tool for constructing transversal Bishop discs, which
may have further applications in the almost complex analysis
and geometry.

This paper was written when the first named author visited the
University of Illinois at Urbana-Champaign during
the Spring semester 2005. He thanks this institution
for hospitality and excellent conditions for work.
In conclusion, the authors thank the referee for many
useful remarks.

\section{Preliminaries}
 In this section we briefly recall some basic properties of almost complex manifolds.

\subsection{Almost complex manifolds}

Let $(M,J)$ be a $C^\infty$-smooth almost complex manifold.
Everywhere below we denote by $\D$  the unit
disc in $\cc$ and by $J_{st}$  the standard complex structure in $\cc^n$;
the value of $n$ is usually clear from the context.
Let $f$ be a smooth map from $\D$ into $M$. We say that $f$ is {\it
 $J$-holomorphic}  if $df \circ J_{st} = J \circ df$. We call such a map $f$
a $J$-{\it holomorphic} disc or a {\it pseudoholomorphic} disc.

The following frequently used statement shows that an  almost complex manifold
$(M,J)$ of complex dimension $n$ can be  locally viewed  as the unit ball $\B$ in
$\cc^n$ equipped with a small almost complex
deformation of $J_{st}$.
\begin{e-lemme}
\label{lemma1}
Let $(M,J)$ be an almost complex manifold of complex dimension $n$. Then for each $p \in
M$,  each  $\delta_0 > 0$, and  each   $k
\geq 0$
 there exist a neighborhood $U$ of $p$ and a
smooth coordinate chart  $Z: U \longrightarrow \B$ such that
$Z(p) = 0$, $dZ(p) \circ J(p) \circ dZ^{-1}(0) = J_{st}$,  and the
direct image $Z_*(J) := dZ \circ J \circ dZ^{-1}$ satisfies
the inequality
$\vert\vert Z_*(J) - J_{st}
\vert\vert_{C^k(\bar {\B})} \leq \delta_0$.
\end{e-lemme}
\proof There exists a diffeomorphism $z$ from a neighborhood $U'$ of
$p \in M$ onto $\B$ such that  $Z(p) = 0$ and $dZ(p) \circ J(p)
\circ dZ^{-1}(0) = J_{st}$. For $\delta > 0$ consider the isotropic dilation
$d_{\delta}: t \mapsto \delta^{-1}t$ in $\cc^n$ and the composite
$Z_{\delta} = d_{\delta} \circ Z$. Then $\lim_{\delta \rightarrow
0} \vert\vert (Z_{\delta})_{*}(J) - J_{st} \vert\vert_{C^k(\bar
{\B})} = 0$. Setting $U = Z^{-1}_{\delta}(\B)$ for  positive
$\delta$  small enough, we obtain the desired result.
\bigskip

Let $(M,J)$ be an almost complex manifold.  Denote by $TM$ the real
tangent bundle of $M$ and by $T_{\cc} M$ its complexification. Recall
that $T_{\cc} M = T^{(1,0)}M \oplus T^{(0,1)}M$ where
$T^{(1,0)}M:=\{ X \in T_{\cc} M : JX=iX\} = \{\zeta -iJ \zeta: \zeta \in TM\}$,
and $T^{(0,1)}M:=\{ X \in T_{\cc} M : JX=-iX\} = \{\zeta + iJ \zeta: \zeta \in TM\}$.
Let $T^*M$ denote the cotangent bundle of  $M$.
Identifying $\cc \otimes T^*M$ with
$T_{\cc}^*M:=\Hom(T_{\cc} M,\cc)$ we consider the set of complex
forms of type $(1,0)$ on $M$:
$
T_{(1,0)}M=\{w \in T_{\cc}^* M : w(X) = 0, \forall X \in T^{(0,1)}M\}
$
and  the set of complex forms of type $(0,1)$ on $M$:
$
T_{(0,1)}M=\{w \in T_{\cc}^* M : w(X) = 0, \forall X \in T^{(1,0)}M\}
$.
Then $T_{\cc}^*M=T_{(1,0)}M \oplus T_{(0,1)}M$.
So we  define the operators $\partial_J$ and
$\bar{\partial}_J$ on the space of smooth functions  on
$M$~: for a  smooth complex function $u$ on $M$ we set
$\partial_J u =  du_{(1,0)} \in T_{(1,0)}M$
and $\bar{\partial}_Ju = du_{(0,1)}
\in T_{(0,1)}M$. As usual, differential forms of any bidegree
$(p,q)$ on $(M,J)$ are defined
by  exterior multiplication.

As usual, an upper semicontinuous function $u$ on $(M,J)$ is called
{\it $J$-plurisubharmonic} on $M$ if the composition $u \circ f$
is subharmonic on $\D$ for every $f \in {\mathcal O}_J(\D,M)$.

Let $u$ be a $\CC^2$ function on $M$, let $p \in M$
and $v \in T_pM$. {\it The Levi
form} of $u$ at $p$ evaluated on $v$ is defined by the equality
$L^J(u)(p)(v):=-d(J^* du)(v,Jv)(p)$.

The following result is well known (see, for instance, \cite{ivro}).
\begin{e-pro}\label{PROP1}
Let $u$ be a $\CC^2$ real valued function on $M$,
let $p \in M$ and $v \in T_pM$.
Then $ L^J(u)(p)(v) = \Delta(u \circ f)(0)$
where $f$ is an arbitrary $J$-holomorphic disc in $M$ such that
 $f(0) = p$ and $
df(0)(\partial / \partial \Re \zeta) = v$ (here $\zeta$ is the
standard complex coordinate variable  in $\cc$).
\end{e-pro}

The Levi form is  invariant with respect to $J$-biholomorphisms.
More precisely, let $u$ be a $\CC^2$ real valued function on $M$,
let $p \in M$ and $v \in T_pM$.
If $\Phi$ is a $(J,J')$-holomorphic diffeomorphism from $(M,J)$ into $(M',J')$,
then $L^J(u)(p)(v) =L^{J'}(u \circ \Phi^{-1})(\Phi(p))(d\Phi(p)(v))$.

Finally, it follows from Proposition~\ref{PROP1} that a
$C^2$-smooth real function $u$ is
$J$-pluri-subharmonic on $M$ if and only if $ L^J(u)(p)(v) \geq 0$
for all   $p \in M$, $v \in T_pM$.
Thus, similarly to  the case of
the integrable structure one arrives in a natural way to
the following definition: a $\CC^2$ real valued
function $u$ on $M$ is {\it strictly
$J$-plurisubharmonic} on $M$ if  $ L^J(u)(p)(v)$
is positive for every $p \in M$, $v \in T_pM \backslash \{0\}$.

It follows easily from Lemma \ref{lemma1} that
for every point $p\in M$ there exists a neighborhood $U$ of $p$
and a diffeomorphism $Z:U \rightarrow \B$ with
center   at $p$ (in the sense that $Z(p) =0$) such that the function
$|Z|^2$ is $J$-plurisubharmonic on $U$ and $Z_*(J) = J_{st} + O(\vert Z \vert)$.

Let $E$ be a real submanifold of codimension $m$ in an almost complex manifold
$(M,J)$ of complex dimension $n$. For every $p$ we denote by
$H_p^J(E)$ the maximal complex (with respect to $J(p)$) subspace
of the tangent space $T_p(E)$. Similarly to the integrable case,
$E$ is said to be a CR manifold if the complex dimension of
$H_p^J(E)$ is independent on $p$; it is called the CR dimension of
$E$ and is denoted by $ \dim_{CR} E$. As usual, by a {\it generic}
submanifold of a complex manifold one means a submanifold $E$ such that
at every point $p \in E$ the complex linear span of $T_p(E)$ coincides with
the tangent space of the ambient manifold.

If $E$ is defined as the common zero set of real functions $\rho_1,...,\rho_m$,
then after the standard identification of $TM$ and $T^{(1,0)}M$, the space
$H_p^J(E)$ can be defined as  the zero subspace of the forms
${\partial}_J\rho_1$, ... ${\partial}_J\rho_m$.
In the present paper we only
deal with the case $m=1$ that is when $E$ is a real hypersurface.

Similarly to the integrable case, by {\it the Levi form} of a real hypersurface
$E = \{ \rho = 0 \}$ at $p \in E$ we mean the conformal class of the Levi form
$L^J(\rho)(p)$ of the defining function $\rho$ on the holomorphic tangent space $H^J_p(E)$.
It is well-known that
\begin{eqnarray}
\label{Leviformula}
L^J_p(\rho)(X_p) = J^*d\rho[X,JX]_p,
\end{eqnarray}
where a vector field $X$ is a smooth section of the $J$-holomorphic tangent bundle $H^J(E)$
of $E$ such that $X(p) = X_p$ for a given vector $X_p \in H^J_p(E)$

\subsection{Normal form of an almost complex structure along a pseudoholomorphic disc}
Throughout the paper, we often use the standard notation
$Z_\zeta:=\frac{\partial Z}{\partial \zeta}$ and
$Z_{\overline\zeta}:=\frac{\partial Z}{\partial \overline \zeta}$.
We don't make a difference between
$\overline Z_{\overline \zeta}$ and $\overline{Z_\zeta}$
as well as between
$\overline Z_\zeta$ and $\overline{Z_{\overline\zeta}}$.

Let $J$ be a smooth almost complex structure on a neighborhood
of the origin in $\cc^2$
and $J(0) = J_{st}$. Denote by
$Z = \left(\begin{array}{cl}
 z\\
 w
\end{array}
\right)
$ the standard coordinates in $\cc^2$.
Then a map $Z:\D \longrightarrow \cc^2$ is $J$-holomorphic if
and only if
it satisfies the following equation
\begin{eqnarray}
\label{basicequation}
Z_{\overline\zeta} - A(Z)\overline{Z_\zeta} = 0,
\end{eqnarray}
where $A(Z)$ is the complex $2\times 2$ matrix  defined by
\begin{eqnarray}
\label{matrixA}
A(Z)v = (J_{st} + J(Z))^{-1}(J_{st} - J(Z))(\overline v).
\end{eqnarray}
It is easy to see that the right-hand side
is $\cc$-linear in $v\in\cc^2$ with respect to the standard
structure $J_{st}$, hence $A(Z)$ is well defined.
Since $J(0) = J_{st}$, we have $A(0) = 0$. However, we will
need a more precise choice of coordinates imposing additional
restrictions on $A$.  We first derive a rule of
transformation of $A$ under diffeomorphisms.

\begin{e-lemme}
\label{transformA}
Let $A(Z)$ be the matrix defined by (\ref{matrixA}).
Let $Z' = Z'(Z)$ be a diffeomorphic change of coordinates
such that
$({\overline Z'}_{\overline Z} + {\overline Z'}_{Z}A)^{-1}$
exists. Then in the new coordinates
\begin{eqnarray*}
A'(Z') = (Z'_ZA + Z'_{\overline Z})
({\overline Z'}_{\overline Z} + {\overline Z'}_{Z}A)^{-1}.
\end{eqnarray*}
\end{e-lemme}
\proof
Consider a $J$-holomorphic disc $\zeta \mapsto Z(\zeta)$
and the disc $Z'(\zeta) = Z'(Z(\zeta))$.
Then
$Z_{\overline\zeta} = A(Z) {\overline Z}_{\overline\zeta}$ and
$Z'_{\overline\zeta}= A'(Z){\overline Z'}_{\overline\zeta}$.
Then we have
$$
(Z'_ZA+ Z'_{\overline Z}){\overline Z}_{\overline\zeta} ={}
Z'_ZZ_{\overline\zeta}+Z'_{\overline Z}\overline Z_{\overline\zeta}={}
Z'_{\overline\zeta} ={}
A'\overline Z'_{\overline\zeta}={}
A'(\overline Z'_{\overline Z}\overline Z_{\overline\zeta}+
\overline Z'_ZZ_{\overline\zeta})={}
A'({\overline Z'}_{\overline Z}+{\overline Z'}_ZA)
\overline Z_{\overline\zeta}.
$$
Since ${\overline Z}_{\overline\zeta}$ is arbitrary,
the desired formula follows.

\begin{e-lemme}
\label{normalization}
Let $Z_0$ be a $J$-holomorphic
disc  close to the disc
$\left(
\begin{array}{cl}
 0\\
 \zeta
\end{array}
\right)$, $\zeta\in\D$.
Then there exists a change coordinates in a neighborhood of
$Z_0(\overline\D)$ such that in the new coordinates we have
$Z_0 = \left(
\begin{array}{cl}
 0\\
 \zeta
\end{array}
\right)$, $\zeta\in\D$.
Moreover, $A(0,\zeta) = 0$,
$A_Z(0,\zeta) = 0$ for $\zeta \in \D$.
\end{e-lemme}

{\bf Remark 1.}
If $A(p)=0$, then the condition $A_Z(p)=0$ means exactly
that for all $J$-holomorphic maps $f:\D\longrightarrow\cc^2$
with $f(0)=p$, in addition to
$f_{\overline\zeta}(0)=0$ we also have $\Delta f(0)=0$.
On the other hand, the integrability condition for $J$
in terms of $A$ is equivalent to certain symmetry
in the expression $A_{\overline Z}+A_Z A$,
therefore for a non-integrable $J$, one generally cannot
achieve $A=0$, $A_{\overline Z}=0$ by a change of coordinates,
even at a point.

{\bf Remark 2.}
If $Z_0\in C^k, 3\le k \le \infty$, then we construct a change
of coordinates of class $C^{k-2}$. We don't think this smoothness
is optimal.

\proof
After a local change of coordinates we have
$Z_0 = \left(
\begin{array}{cl}
 0\\
 \zeta
\end{array}
\right)$.
The $J$-holomorphicity condition of $Z$ implies that in
these coordinates we have
$$(A\circ Z)(\zeta) = \left(
\begin{array}{cll}
 \alpha(\zeta) & 0\\
 \beta(\zeta) & 0
\end{array}
\right).
$$
Consider a local change of coordinates of the form
\begin{eqnarray*}
& &z' = a_{10}z + a_{01} \overline z +   a_{11}z\overline z,\\
& &w' = w + b_{10}z + b_{01} \overline z  + b_{11}z\overline z,
\end{eqnarray*}
where $a_{jk}$, $b_{jk}$ are smooth functions of $w$
with $|a_{10}|\ne|a_{01}|$.
In the new coordinates the disc $Z_0$ does not change.
We have
$$Z'_Z(0,w) = \left(
\begin{array}{cll}
 a_{10} & 0\\
 b_{10} & 1
\end{array}
\right), \qquad
Z'_{\overline Z}(0,w) = \left(
\begin{array}{cll}
 a_{01} & 0\\
 b_{01} & 0
\end{array}
\right).
$$
By Lemma \ref{transformA} we have
\begin{eqnarray*}
A'(Z') = (Z'_ZA + Z'_{\overline Z})
({\overline Z'}_{\overline Z} + {\overline Z'}_{Z}A)^{-1}.
\end{eqnarray*}
The condition $A'(0,w') = 0$ implies
\begin{eqnarray*}
\left(
\begin{array}{cll}
 a_{10} & 0\\
 b_{10} & 1
\end{array}
\right)
\left(
\begin{array}{cll}
 \alpha & 0\\
 \beta & 0
\end{array}
\right) +
\left(
\begin{array}{cll}
 a_{01} & 0\\
 b_{01} & 0
\end{array}
\right) = 0
\end{eqnarray*}
and therefore
\begin{eqnarray*}
& &\alpha a_{10} + a_{01} =0,\\
& &\alpha b_{10} + \beta + b_{01} = 0.
\end{eqnarray*}
Hence the condition $A'(0,w') = 0$ determines the
functions $a_{01}$ and $b_{01}$ for given
$\alpha$, $\beta$, $a_{10}$ and $b_{10}$.

Thus without loss of generality we  assume that
$A(0,w) =0$ and  $a_{01} = b_{01} = 0$.
This also implies that $A_w(0,w) = A_{\overline w}(0,w) = 0$.
So  we need to achieve
$A'_{z'}(0,w') = 0$.

We have $A'_z = A'_{z'}z'_z + A'_{\overline z'}{\overline z'}_z$.
Since $a_{01} = b_{01} = 0$,
then $Z'_{\overline Z}(0,w) = 0$.
Then $A'_z(0,w) = a_{10}(w)A'_{z'}(0,w')$.
Hence $A'_{z'}(0,w') = 0$ if and only if $A'_z(0,w) = 0$.
Since $A(0,w) = 0$, then for $z = 0$ we have
\begin{eqnarray*}
A'_z(0,w) = (Z'_ZA_z + Z'_{\overline Z z})({\overline Z'}_{\overline Z} + {\overline Z'}_{Z}A)^{-1}.
\end{eqnarray*}
We also have
\begin{eqnarray*}
Z'_{\overline Z z}(0,w) = \left(
\begin{array}{cll}
 z'_{\overline z z} & z'_{\overline w z}\\
 w'_{\overline z z} & w'_{\overline w z}
\end{array}
\right) = \left(
\begin{array}{cll}
 a_{11} & (a_{10})_{\overline w}\\
 b_{11} & (b_{10})_{\overline w}
\end{array}
\right).
\end{eqnarray*}
Put
\begin{eqnarray*}
A_z(0,w) = \left(
\begin{array}{cll}
 \alpha & \beta\\
 \gamma & \delta
\end{array}
\right).
\end{eqnarray*}
Then the condition $A'_{z'}(0,w') = 0$ takes the form

\begin{eqnarray*}
\left(
\begin{array}{cll}
 a_{10} & 0\\
 b_{10} & 1
\end{array}
\right)
\left(
\begin{array}{cll}
 \alpha & \beta\\
 \gamma & \delta
\end{array}
\right) +
\left(
\begin{array}{cll}
 a_{11} & (a_{10})_{\overline w}\\
 b_{11} & (b_{10})_{\overline w}
\end{array}
\right) = 0,
\end{eqnarray*}
so the coefficients $a_{11}$, $b_{11}$ are determined by $a_{10}, b_{10}, \alpha,\gamma$.
The functions $a_{10}$ and $b_{10}$ are found as solutions
of the following classical elliptic equations (see \cite{Ve}):
\begin{eqnarray*}
& &(a_{10})_{\overline w} + \beta a_{10} = 0,\\
& &(b_{10})_{\overline w} + \beta b_{10} + \delta = 0.
\end{eqnarray*}
This completes the proof of the lemma.
\bigskip

Let $E$ be a real hypersurface through the origin in $\cc^2$
equipped with a smooth almost complex structure
$J$. Even if $J(0) = J_{st}$, the Levi form of $E$ with respect
to $J$ at the origin does not necessarily coincide with
the Levi form with respect to $J_{st}$.
However, if the coordinates are normalized according
to the previous lemma, then the Levi forms with respect
to $J$ and $J_{st}$ are the same.

\begin{e-lemme}
\label{Leviform}
Assume that $A(0) = A_Z(0) = 0$. Then
the Levi form of $E$ at the origin with respect to the structure $J$ coincides with the Levi form of $E$
at the origin with respect to the structure $J_{st}$.
\end{e-lemme}
\proof We can assume that $E$ is given by
\begin{eqnarray*}
\rho(Z) = \Re z + Q(Z) + H(Z) + o(\vert Z \vert^2) = 0,
\end{eqnarray*}
where $Q(Z)$ is a complex quadratic form and $H(Z)$ is a hermitian quadratic form. It is sufficient to consider
the Levi form of $\rho$ at the origin evaluated on
the vector $V = (0,1)$.
So we consider a $J$-holomorphic disc $Z(\zeta)$ satisfying
$Z(\zeta)=  \left(
\begin{array}{cl}
 z(\zeta)\\
 w(\zeta)
\end{array}
\right) = \left(
\begin{array}{cl}
 0\\
 \zeta
\end{array}
\right) + O(\vert \zeta \vert^2)$
(the existence of such a disc follows by the classical
Nijenhuis-Woolf theorem, see for instance  \cite{Si}).
It follows from the assumption on the matrix $A$ and
the $J$-holomorphicity
of $Z$ that $w_{\overline \zeta} = a\overline \zeta
+ o(\vert \zeta \vert)$ so that
$w_{\zeta\overline\zeta}(0) = 0$.
Therefore $\Delta (\rho \circ Z)(0) = H(V)$
which proves the lemma.

\section{Bishop discs and the Bishop equation}
Let $(M,J)$ be a
smooth almost complex
manifold of real dimension
$2n$ and $E$  a generating submanifold of $M$ of real codimension
$m$. A $J$-holomorphic disc $f:\D \longrightarrow M$ continuous on
$\overline \D$ is called a {\it Bishop disc} if $f(b\D) \subset E$,
where $b\D$ denotes the boundary of $\D$.
The existence and local parametrization of  certain classes of Bishop discs
attached to $E$ are obtained in \cite{KS}. Here we give a more precise description of
small Bishop discs which will be used in our constructions.

\subsection{Bishop's equation as the Riemann-Hilbert type problem for elliptic PDE systems}
Let $E$ be a smooth generic submanifold in a smooth (always supposed $C^{\infty}$)
almost complex manifold $(M,J)$  defined as the zero set of an
$\R^m$-valued  function $\rho = (\rho^1,...,\rho^m)$  on $M$.
Then a smooth map
$$f:\D \longrightarrow M$$
$$f:\zeta \mapsto f(\zeta)$$  continuous on
$\overline\D$ is a Bishop disc if and only if it satisfies the
following non-linear boundary problem of the
Riemann-Hilbert type for the
quasi-linear operator $\overline\partial_J$:

\begin{displaymath}
(RH): \left\{ \begin{array}{ll}
\overline\partial_J f(\zeta) = 0, \,\,\,\zeta \in \D\\
\rho(f)(\zeta) = 0, \,\,\, \zeta \in b\D
\end{array} \right.
\end{displaymath}

In order to obtain a local  description of  solutions of this problem we
fix a chart $U \subset M$ and a coordinate diffeomorphism
$Z: U \longrightarrow \B$ where $\B$ is  the unit ball
of $\cc^n$. Identifying $M$ with $\B$
we can assume that in these coordinates $J= J_{st} + O(\vert Z \vert)$
and the norm $\parallel J - J_{st} \parallel_{C^k(\overline \B)}$ is
small enough for some  positive real
$k$ in accordance with Lemma \ref{lemma1}. (Here $k > 1$ can be chosen arbitrary;
 we  assume it for convenience to be real positive {\it non-integral}
and fix it throughout what follows.)
 More precisely,
using the notation $Z = (z,w)$, $z = (z_1,...,z_m)$, $z = x + iy$,
$w = (w_{1},...,w_{n-m})$ for the standard coordinates in $\cc^n$, we can also
assume that $E \cap U$ is described by the equations

\begin{eqnarray}
\label{manifold}
\rho(Z) = x - h(y,w) = 0
\end{eqnarray}
with  vector-valued $C^{\infty}$-function $h:\B \longrightarrow \R^m$ such that
$h(0) = 0$ and $\bigtriangledown h(0) = 0$.

Similarly to  the proof of Lemma \ref{lemma1} consider the isotropic dilations
$d_{\delta}: Z \mapsto Z' = \delta^{-1}Z$. In the new
$Z$-variables  (we drop the primes) the image
$E_{\delta} = d_{\delta}(E)$ is
defined by the equation $\rho_\delta (Z):\delta^{-1}\rho(\delta Z) = 0$. Since the function
$\rho_\delta$ approaches $x$ as $\delta \longrightarrow 0$,
the manifolds $E_\delta$ approach the flat manifold $E_0=\{
x = 0\}$, which, of course, may be identified with the real tangent
space to $E$ at the origin. Furthermore, as         seen in the proof
of Lemma \ref{lemma1}, the structures $J_\delta :(d_\delta)_*(J)$
converge to $J_{st}$ in the $C^k$-norm  as
$\delta \longrightarrow 0$. This allows us to find explicitly   the
$\overline\partial_J$-operator in the $Z$ variables.

Consider now a $J_\delta$-holomorphic disc
$$Z: \D \longrightarrow (\B,J_\delta)$$
$$Z: \zeta \mapsto Z(\zeta)$$
 of class $C^k(\overline D)$. The
$J_\delta$-holomorphicity condition $J_\delta(Z) \circ dZ = dZ \circ J_{st}$
can be written in the following form.
\begin{eqnarray}
\label{CR}
Z_{\overline\zeta} - A_{J,\delta}(Z)
{\overline Z}_{\overline \zeta}
= 0,
\end{eqnarray}
where $A_{J,\delta(Z)}$ is the
complex $n\times n$ matrix  of an operator the composite of which with complex
conjugation is equal to the endomorphism
$(J_{st} + J_\delta(Z))^{-1}(J_{st} - J_\delta(Z))$
(which is an anti-linear operator with respect to the standard
structure $J_{st}$). Hence the entries of the matrix $A_{J,\delta}(Z)$
are smooth functions of $\delta,Z$ vanishing identically in $Z$
for $\delta = 0$.

Recall that the Cauchy-Green transform is defined by
\begin{eqnarray}
\label{CauchyGreen}
Tf(\zeta) = \frac{1}{2\pi i} \int\int_{\D} \frac{f(\tau)}{\tau -
  \zeta}d\tau \wedge d\overline\tau.
\end{eqnarray}
It is well-known  (\cite{Ve}, p. 56, Theorem 1.32) that
$T:C^k(\overline \D) \longrightarrow C^{k+1}(\overline \D)$
is a bounded operator when $k$ is non-integer.
Furthermore, $(Tf)_{\overline \zeta} = f$.
Therefore we can write  $\overline\partial_J$-equation (\ref{CR}) as follows:
$$
[ Z - T ( A_{J,\delta}(Z)
{\overline Z}_{ \overline \zeta}  )  ]_{\overline\zeta} = 0.
$$
The arising non-linear operator
$$
\Phi_{J,\delta}: Z \longrightarrow \tilde Z =  Z - T [ A_{J,\delta}(Z)
{\overline Z}_{ \overline \zeta}]
$$
takes the space   $C^{k}(\overline\D)$  into itself.
Thus, $Z$ is $J_\delta$-holomorphic disc if
and only if $\Phi_{J,\delta}(Z)$ is a holomorphic disc (in the usual sense)
on $\D$.
For sufficiently small positive $\delta$  the map $\Phi_{J,\delta}$ is an
invertible operator on a neighborhood  of zero in $C^k(\overline D)$
which  establishes a one-to-one correspondence between the
sets of $J_\delta$-holomorphic and holomorphic discs  in $\B^n$.

These considerations allow us to replace the non-linear
Riemann-Hilbert problem
(RH) by {\it the generalized Bishop equation}
\begin{eqnarray}
\label{Bishop2}
\rho_\delta(\Phi_{J,\delta}^{-1}(\tilde Z))(\zeta) = 0, \;
\zeta \in b\D,
\end{eqnarray}
for  an unknown {\it holomorphic} function $\tilde Z$ in $\D$.

If $\tilde Z$ is a solution of the boundary  problem (\ref{Bishop2}),
then $Z = \Phi_{J,\delta}^{-1}(\tilde Z)$ is a Bishop disc with  boundary
attached to $E_\delta$. Since the manifold $E_\delta$ is
biholomorhic via  isotropic dilations to the initial manifold $E$,
the solutions of the equation (\ref{Bishop2}) allow to describe
the Bishop discs attached to $E$. Of course, this gives just the discs
{\it close enough} in the $C^k$-norm to the trivial solution $Z \equiv 0$
of the problem (RH). We will call such discs {\it small}.

\subsection{Schwarz-Green formula and Bishop's equation}
This is also useful to give a more explicit form of the Bishop equation as a non-linear
system of singular integral equations.
Let

\begin{eqnarray*}
P_0f = \frac{1}{2 \pi i} \int_{b\D} f(\tau) \frac{d\tau}{\tau}
\end{eqnarray*}
denote the average value of $f$ on $b\D$ and let
\begin{eqnarray*}
Su(\zeta) = \frac{1}{2\pi i}\int_{b\D} \frac{\tau +
  \zeta}{\tau - \zeta} u(\tau)\frac{d\tau}{\tau}
\end{eqnarray*}
denote the Schwarz integral. Consider also the Cauchy integral

\begin{eqnarray}
\label{Cauchy}
Kf(\zeta) = \frac{1}{2 \pi i}\int_{b\D} \frac{f(\tau)d\tau}{\tau - \zeta}.
\end{eqnarray}
We note
$$S = 2K - P_0.$$

We begin with a version of the Cauchy-Green formula replacing the Cauchy kernel by the Schwarz kernel;
this is a variation of the classical results \cite{Ve}.

\begin{e-pro}
Let $f = u + iv$ a complex function of class $C^1$ on the unit disc $\D$ of $\cc$. Then for every $\zeta \in \D$
we have
\begin{eqnarray*}
f(\zeta) = Su(\zeta) + iv_0 +  T f_{\overline\zeta}(\zeta) - \overline{T f_{\overline\zeta}(1/\overline\zeta)},
\end{eqnarray*}
where $v_0 = P_0v$.
\end{e-pro}
\proof Let
\begin{eqnarray*}
g(\zeta) = f(\zeta) - Su(\zeta) - iv_0 - T f_{\overline\zeta}(\zeta) +
\overline{T f_{\overline\zeta}(1/\overline\zeta)}.
\end{eqnarray*}
Then the function $g$ is holomorphic on $\D$ because
$\partial_{\overline\zeta}
T f_{\overline\zeta}(\zeta)=f_{\overline\zeta}(\zeta)$
and the functions
$Su$ and $\overline{T f_{\overline\zeta}(1/\overline\zeta)}$
are holomorphic. Furthermore,
$\Re g \vert_{b\D} = 0$ since $\Re Su\vert_{b\D} = u$
(in the sense of limiting boundary values from inside) and
$1/\overline\zeta = \zeta$,
$Tf_{\overline\zeta}(\zeta)=Tf_{\overline\zeta}(1/\overline\zeta)$
for $\zeta \in b\D$. Hence $g = ic$ where $c$ is a real constant.

In order to prove that $c = v_0$,
it suffices to show that $P_0T f_{\overline\zeta} = 0$.
Note that for for every bounded $g$, in particular for
$g=f_{\overline\zeta}$ the function $Tg$ is
holomorphic in $\cc\setminus\D$ and vanishes at the infinity.
Hence $P_0Tg = 0$, and the proposition follows.
\bigskip

Let now $E$ be a generic submanifold of codimension $n - k$
in an almost complex manifold $(M,J)$ of complex dimension $n$.
Fixing local coordinates, we can assume that
$E$ is a submanifold of $\cc^n$ through the origin, $J(0) = J_{st}$.
Similarly to \cite{Tu2} we can also assume that $E$ is given by the following parametric
equations
\begin{eqnarray*}
\Re Z = h(\Im Z,t),
\end{eqnarray*}
where $Z \in \cc^n$ are the standard complex coordinates and $t \in \R^k$ is a parameter.
Furthermore,
\begin{eqnarray*}
h(0,0) = 0,\;\;
\frac{\partial h}{\partial \Im Z}(0,0) = 0,\;\;
\rank \frac{\partial h}{\partial t}(0,0) = k.
\end{eqnarray*}

Let $Z: \zeta \mapsto Z(\zeta)$ be a $J$-holomorphic Bishop disc for $E$ in a sufficiently small neighborhood
of the origin.
Then it satisfies the $J$-holomorphicity equations
\begin{eqnarray*}
Z_{\overline\zeta} - A(Z(\zeta))\overline{Z_{\zeta}(\zeta)} = 0.
\end{eqnarray*}
Since the disc $Z$ takes its values in a neighborhood of the origin small enough, the norm of the matrix $A$ also
is supposed to be small. The boundary condition $Z(b\D) \subset E$ means that
\begin{eqnarray*}
\Re Z(\zeta) = h(\Im Z(\zeta),t(\zeta))
\end{eqnarray*}
for some function $t(\zeta)$ on $b\D$. Set
\begin{eqnarray*}
Pf(\zeta) = Tf(\zeta) - \overline{Tf(1/\overline\zeta)}.
\end{eqnarray*}
Then we obtain that the above conditions are equivalent to   the following Bishop equation for the map $Z$:
\begin{eqnarray}
\label{Schwarz-Bishop}
Z = S h(\Im Z, t) + ic_0 + P(A(Z)\overline{Z_{\zeta}}).
\end{eqnarray}
For a given non-integral $\alpha > 1$, sufficiently small function $t(\zeta) \in (C^\alpha(b\D))^k$,
and $c_0 \in \R^n$
this equation has a unique solution of class
$C^\alpha(\overline \D)$ by the implicit function theorem.

In particular, if $E$ is given by the equations

\begin{eqnarray*}
& &x = \varphi(y,w),\\
& &z = x + i y \in \cc^{n-k}, \; w = u +i v \in \cc^k,\\
& &\varphi(0,0) =0, \; d\varphi(0,0) = 0,
\end{eqnarray*}
then

$$
\left(
\begin{array}{cl}
x\\
u
\end{array}
\right) = h(y,v,t) = \left(
\begin{array}{cl}
\varphi(y,t + i v)\\
t
\end{array}
\right)
$$
can be used for Bishop's equation  (\ref{Schwarz-Bishop}).

We can modify Bishop's equation to define
$Z(1) = Z_0 = X_0 + iY_0$, $X_0 = h(Y_0,t_0)$.
Put
\begin{eqnarray*}
& & S_1u = Su - Su(1),\\
& &P_1f = Pf - (Pf)(1).
\end{eqnarray*}
Then we have the modified Bishop's equation
\begin{eqnarray}
\label{modifBishop}
Z = S_1h(\Im Z, t) + Z_0 + P_1(A \overline{Z_{\zeta}}),
\end{eqnarray}
where $t(1) = t_0$. The solution satisfies $Z(1) = Z_0$.

We point out that it follows immediately from the equation (\ref{Schwarz-Bishop}) or (\ref{modifBishop})
that Bishop's discs depend smoothly on
deformations of the almost complex structure $J$.

\subsection{Parametrization of  discs}

Let $p $ be a point of a generic submanifold $E$ in an almost complex manifold $(M,J)$.
Fix local coordinates so that $E$ is defined by (\ref{manifold})
and $p = 0$.
Our goal is to describe the solutions $\tilde Z$ of the generalized
Bishop equation (\ref{Bishop2}) satisfying the condition
$\Phi_{J,\delta}^{-1}(\tilde Z)(1) = 0$.
Our argument is similar to \cite{KS}.

Let  ${\cal U}$ be a  neighborhood of the origin in $\R$,
$X'$ a sufficiently small neighborhood of the origin in the Banach  space
$(C^k(\overline\D) \cap {\cal O}(\D))^m$ (with positive non-integral $k > 1$),
$X''$ a neighborhood of the origin in the Banach space
$(C^k(\overline\D) \cap {\cal O}(\D))^{n-m}$, and $Y$ the Banach space
$(C^{k}(b\D))^m$.
If $\tilde z \in X'$, $\tilde z:\zeta \mapsto \tilde z(\zeta)$ and
$\tilde w \in X''$, $\tilde w:\zeta \mapsto \tilde w(\zeta)$ are holomorphic discs,
then we denote by $\tilde Z$ the holomorphic disc $\tilde Z = (\tilde z,\tilde w)$.

Set $X = X' \times X''$.
Given $\tilde Z \in X \subset (C^k(\overline\D) \cap {\cal O}(\D))^n$,
we put $(z_\delta,w_\delta) =  \Phi_{J,\delta}^{-1}(\tilde Z)$ and consider the map of Banach spaces
$R:  X \times {\cal U}  \longrightarrow Y \times \R^m \times \cc^{n-m}$
defined as follows:
$$R : (\tilde z,\tilde w,\delta) \mapsto
\left ( \rho_\delta(\Phi_{J,\delta}^{-1}(\tilde Z))(\bullet)\vert_{ b\D}, y_\delta(1),w_\delta(1) \right ).$$

Let $\phi$ be a   $C^{2k}$-map between two domains in $\R^n$ and $\R^m$ ; it determines a map
$\omega_\phi$ acting   by composition on  $C^k$-smooth maps $g$ into
the source domain:
$\omega_\phi: g \mapsto \phi(g)$. The well-known fact is that
$\omega_\phi$ is a $C^k$-smooth map between the corresponding spaces
of $C^k$-maps. In our case this means that the map $R$ is of class $C^k$.
\begin{e-lemme}
The tangent map $$DR:= D_{X}R(0,0,0): X
\longrightarrow Y \times \R^m \times \cc^{n-m}$$ (the partial derivative with respect
to the space $X$) is surjective.
\end{e-lemme}
\proof This is easy to see that the  map $DR$ defined  by the equality
$$DR(\hat Z) = (\Re \hat z_1,...,\Re \hat z_m,   \hat y(1), \hat w(1))$$
for any $\hat Z = (\hat z,\hat w) \in X$.
Recall that the Hilbert transform $H$ can be defined as a singular integral operator
\begin{eqnarray*}
Hu(\zeta) = \frac{1}{2\pi i} (p.v.) \int_{0}^{2\pi}
\frac{e^{i\theta} + \zeta}{e^{i\theta}-\zeta}
u(e^{i\theta})d\theta.
\end{eqnarray*}
The operator  $H$ is a bounded linear operator on $C^k(b\D)$
for any non-integral positive $k$.
Set  $H_1u = Hu - (Hu)(1)$ so
that $(H_1u)(1) = 0$.
Now for given $u \in (C^k(b\D))^m$, given $ c \in \R^m$ and
$q \in \cc^{n-m}$
we set $\hat z = u + iH_1u + ic$ and denote by $\hat w$
a holomorphic function satisfying
$\hat w(1) = q$.
Then $DR( \hat z, \hat w) = (u,c,q)$ so the map $DR$ is surjective. Moreover, its kernel is
canonically isomorphic
to the space $(C^k(\overline\D) \cap {\mathcal O}_1(\D))^{n-m}$ of holomorphic functions $\hat w$ of class
$C^k(\overline \D)$ satisfying $\hat w(1) = 0$.
\bigskip

Now given $\hat w$ we look for a solution $\tilde Z$ of the equation (\ref{Bishop2}). By the
implicit function theorem there exist
$\delta_0 > 0$,
a neighborhood $V_1$ of the origin in $X'$,
a neighborhood $V_2$ of the origin in $(C^k(\overline\D) \cap {\mathcal O }_1(\D))^{n-m}$,
 and
a $C^k$ smooth  map
$G:  V_2  \times [0,\delta_0] \longrightarrow V_1$
such that for every $(\hat w,\delta) \in V_2  \times
[0,\delta_0]$ the holomorphic function $\tilde Z = G(\hat w,\delta)(\bullet)$ is
the unique solution of generalized Bishop's equation (\ref{Bishop2})
belonging to $V_1 \times V_2$ and satisfying the condition
$\Phi_{J,\delta}^{-1}(\tilde Z)(1) = 0$.

Now, the pullback $Z = \Phi_{J,\delta}^{-1}(\tilde Z)$ gives us
a $J_\delta$-holomorphic disc attached
to $E_\delta$ and satisfying $Z(1) = 0$.   Thus, the  initial data consisting of
 $\hat w \in (C^k(\overline\D) \cap {\mathcal O}_1(\D))^{n-m}$ define
for each small $\delta$ a unique
$J_\delta$-holomorphic disc $Z_0 = Z(\hat w)(\bullet)$
attached to $E_\delta$ and satisfying $Z_0(1) = 0$. We will call this disc $Z_0$ {\it the lift} of $\hat w$
attached at $p = 0$ (we write $Z_p(\hat w)$ for an arbitrary $p$).
 Since the almost complex structures $J$
and $J_\delta$ are biholomorphic via  isotropic dilations, we can
give   the following description of local solutions of  Bishop's
equation (the problem (RH)).

\begin{e-pro}
\label{theosolvingBishop}  The set ${\mathcal A}^J_p$ of
$J$-holomorphic Bishop discs satisfying
$Z(1) = p$
is a Banach manifold and its tangent space at $Z  \equiv 0$ is
canonically isomorphic to $(C^k(\overline\D) \cap
{\mathcal O}_1(\D))^{n-m}$.
Moreover, the parametrization map
$(C^k(\overline\D) \cap {\mathcal O}_1(\D))^{n-m} \times
E \ni(\hat w,p)\mapsto Z_p[\hat w]$ is smooth.
\end{e-pro}

The proof follows from the above analysis of the Bishop equation. One
merely fixes    some value of $\delta$, $0<\delta\leq \delta_0$
and observes that the families
of  Bishop discs corresponding to distinct values of
$\delta\neq 0$ are taken into one another by the corresponding dilations.

In what follows we say that a $J$-holomorphic Bishop disc $Z$ of $E$
satisfying $Z(1) = p$ is {\it attached at $p \in E$}.

The representation (\ref{modifBishop}) of Bishop's equation is convenient since the structure $J$ and
defining functions of $E$
appear explicitly in the coefficients of integral equations. This implies immediately that solutions
depends smoothly on deformations of the structure. In particular,   if $E_1$ and $E_2$
are $C^{2k}$-close submanifolds of $\cc^n$ defined by equations of the form
(\ref{manifold}) and $J_1$ and $J_2$ are $C^{2k}$
close almost complex structures, then  there exists a (locally defined)
diffeomorphism between the corresponding    $C^{k}$ Banach submanifolds
of Bishop discs ${\mathcal A}^{J_1}_0$ and ${\mathcal A}^{J_2}_0$
that depends smoothly on the pairs $E_j$, $J_j$, $j=1,2$,
and is the identity  in the case of equal pairs.

Moreover, if we apply the implicit function theorem to the equation (\ref{modifBishop}) we obtain that
every Bishop discs is uniquely determined by
the initial data $Z_0$ and $t(\zeta)$. Thus, we obtain another parametrization of the manifold ${\cal A}^J_p$.

\section{Filling by pseudoholomorphic Bishop discs: pseudoconvex and finite type cases}

In this section we consider the case where $E$ is a real hypersurface. The following simple statement will
be crucially used.

\begin{e-pro}
\label{filling}
Suppose that there exists $Z \in {\cal A}_p^J$ such that the normal derivative vector
$dZ(\partial_{\Re \zeta}\vert_1)$
is not tangent to $E$ at $p$ (that is the Bishop disc $Z$ is attached to $E$ at $p$ transversally).
Then Bishop discs from $\cup_{q \in E} {\mathcal A}_q^J$ of $E$ fill a one-sided neighborhood of $p \in E$.
\end{e-pro}

\proof
We can assume that $Z$ is the lift of $\hat w \in (C^k(\overline\D) \cap {\mathcal O}_1(\D))^{n-m}$
(defined in the proof of Proposition
\ref{theosolvingBishop}). It follows by Proposition \ref{theosolvingBishop} that for every point
$q \in E$ in a neighborhood of
$p$ there exists a disc $Z_{q}(\hat w)$ (the lift of $\hat w$ attached to $E$ at $q$) which is
a small deformation of $Z$.
So the normals $Z_q([0,1])$
of these discs at $1$ fill a one-sided neighborhood of $p$. This completes the proof.

\bigskip

{\bf Remark.} If $Z$ is a Bishop disc attached to $E$ at $p = 0$, the tangent space $T_Z({\mathcal A}_p^J)$ of the manifold
of Bishop discs attached at $p$ consists of maps $\dot Z: \overline\D \longrightarrow \cc^n$
with $\dot Z(1) = 0$  satisfying a Riemann-Hilbert type system of the form
\begin{displaymath}
\left\{ \begin{array}{ll}
\dot Z_{\overline\zeta} + a(\zeta) \overline{\dot Z}_{\overline\zeta} + b(\zeta) \dot Z +
c(\zeta) \overline{\dot Z}=0,\;\;
\zeta \in \D,\\
\Re (\partial \rho(Z) \dot Z) = 0,\;\;
\zeta \in b\D,
\end{array} \right.
\end{displaymath}
where the matrix valued coefficients $a$, $b$, $c$ involve
$\overline Z_{\overline\zeta}$ and the values of the matrix $A$ and its derivatives
$A_Z$ and $A_{\overline Z}$ on the disc $Z$.
This system is obtained by varying the Bishop equation.
The solutions $\dot Z$ of this problem are called {\it the infinitesimal variations}
(or {\it infinitesimal perturbations}) of the disc $Z$. In general the study of the above boundary problem is quite complicated. Fortunately, in the present paper we deal with small Bishop discs only.
This implies that the above Bishop disc $Z$ can be chosen small enough. Therefore, the above Riemann-Hilbert problem can be viewed as an arbitrarily  small perturbation of the   usual $\overline\partial$-equation on the unit disc with  the boundary condition
$\Re (\partial \rho(0) \dot Z) = 0$.
In particular the corresponding linear operator is surjective in  suitable functional spaces and the above linearized equation indeed determines the tangent space to the manifold ${\cal A}_p^J$. A neighborhood of the disc $Z$ in ${\cal A}_p^J$ is in a one-to-one correspondence  with a neighborhood of the origin in this tangent space by the implicit function theorem. The last proposition means that if the space of discs $\dot Z$ satisfying the above conditions contains a disc
with the normal derivative vector at 1 transversal to $E$, then Bishop disc fill one-sided neighborhood of $E$.
In this section  we  study the above linear Riemann-Hilbert problem by geometric tools: using suitable
non-isotropic dilations of coordinates we represent it as a small deformation of the corresponding Riemann-Hilbert
problem for the usual $\overline\partial$-operator. An analytic study of this problem is postponed to section 5.

We will need the following statement.

\begin{e-pro}
\label{perturbation}
Let $E$ be a real hypersurface  in an almost complex manifold $M$ with an integrable structure $J_0$.
Assume that $E$ contains no complex hypersurfaces. Then
for any almost complex structure $J$ close enough to $J_0$ in the $C^k$, $k>2$, norm
the $J$-holomorphic Bishop discs
of $E$ fill a one-sided neighborhood of every point of  $E$.
\end{e-pro}
\proof
According to \cite{Tu} the manifold ${\mathcal A}_p^{J_{0}}$ contains a disc transversal to $E$ at $p \in E$.
Since the manifolds ${\mathcal A}_p^{J}$ depend smoothly on $J$, they also contain a transversal disc
if $\parallel J - J_{0} \parallel_{C^k}$ is small enough.

Of course, a similar statement remains true if $E$ is a generic submanifold in $(M,J_0)$ minimal
in the sense of \cite{Tu}.

\subsection{Pseudoconvex hypersurfaces}

As a first consequence we obtain the following

\begin{e-pro}
Let $\rho$ be a $J$-plurisubharmonic function of class $C^2$ on an almost complex manifold $(M,J)$ of complex dimension 2
and the hypersurface $E = \rho^{-1}(0)$ contains no $J$-holomorphic discs (we assume $d\rho \neq 0$ on $E$).
Then there exists a neighborhood $U$ of $E$ such that the Bishop discs
of $E$ fill  $U \cap \{ \rho < 0 \}$.
\end{e-pro}
\proof Consider a (sufficiently small) Bishop disc $Z \in {\cal A}_p^J$. Then
$\rho \circ Z \leq 0$ on $\D$ by the maximum principle. But then $dZ(\partial_{\Re \zeta}\vert_1)$
is transversal to $E$ by the Hopf lemma. So we apply Proposition \ref{filling}.

\bigskip

The above result remains true for arbitrary dimension (see \cite{BaTrRo} for the case of $\cc^n$).

\begin{e-pro}
Let $E$ be a real hypersurface defined as the zero set of a  $J$-plurisubharmonic function $\rho$ of class $C^2$
in an almost complex manifold $(M,J)$.
Suppose that there exists a complex tangent vector field $X$ on $E$ such that
for any $p \in E$  there are no (germs of) $J$-holomorphic discs in $E$ tangent to $X(p)$.
Then there exists a neighborhood $U$ of $E$ such that Bishop discs of $E$ fill $U \cap \{ \rho < 0 \}$.
\end{e-pro}
\proof If the Banach manifold ${\mathcal A}_p^J$ contains no transversal discs,
by the Hopf lemma  all the discs from ${\mathcal A}_p^J$ are
contained in $E$. Similarly to the previous section, fix local coordinates centered at $p$ and consider
isotropic dilations $\Lambda_\delta$.  Fix a point $q \in E_0$ close enough to $p$.
Then the tangent vectors at the centers of Bishop discs (with respect to $J_{st}$ ) attached to $E_0$
at $0$  fill a real sphere in the holomorphic tangent space $H_q^{J_{st}}(E_0)$.
By continuity, a similar property holds for $\delta$ small enough so there is a point in $E_\delta$
admitting a germ of a $J_\delta$-holomorphic disc in any complex tangent direction and contained in $E_\delta$.
This contradiction shows that the manifold ${\mathcal A}_p^{J}$ contains a transversal disc and
we apply Proposition  \ref{filling}.
\bigskip

As a corollary we obtain the following global version of this statement.

\begin{e-cor}
Let $\rho$ be a $C^2$ plurisubharmonic function on an almost  complex manifold.
Suppose that for some $c \in \R$ its level set $E:= \rho^{-1}(c)$ is a compact
hypersurface of class $C^2$ which contains no germs of holomorphic discs.
Then there exists a neighborhood $U$ of $E$ such that Bishop discs of $E$ fill
$U \cap \{ \rho < c \}$.
\end{e-cor}

\subsection{Finite type  hypersurfaces}

Consider a real hypersurface $E$ in an almost complex manifold $M$ of complex dimension $2$.
By its type we mean the
supremum of tangency order of $E$ with regular $J$-holomorphic discs at $p \in E$ (the properties of finite type
hypersurfaces in almost complex manifolds are studied in the recent work   \cite{BaMa}).

\begin{e-pro}
If $E$ is of finite type, then its Bishop discs fill a one-sided neighborhood of every point of $E$.
\end{e-pro}

Let an integral  $m > 1$ be the type of $E$ at the origin. There exists a regular disc $Z$
such that
\begin{eqnarray}
\label{tangency}
(\rho \circ Z)(\zeta) = O(\vert \zeta \vert^m).
\end{eqnarray}
We can choose local coordinates so that
$Z(\zeta) = (0,\zeta)$
and the $J$-holomorphicity conditions for a map $(w,z):\D \longrightarrow \cc^2$ have the form
\begin{eqnarray*}
& &w_{\overline\zeta} + a(w,z)\overline{w}_{\overline\zeta} = 0,\\
& &z_{\overline\zeta} + b(w,z)\overline{z}_{\overline\zeta} = 0,
\end{eqnarray*}
where
\begin{eqnarray}
\label{vanishing}
a(\bullet,0) = 0, \; b(0,0) = 0.
\end{eqnarray}

Indeed, we can identify $M$ with $\cc^2$ equipped with an almost complex structure
$J$ such that $J(0) = J_{st}$. Using the classical result of Nijenhuis-Woolf (see, for instance, \cite{Si})
we construct a foliation of a neighborhood of the origin in $\cc^2$ by a family of $J$-holomorphic discs containing $Z$.
Similarly, consider also a transversal foliation. After a local diffeomorphism the discs of these foliations become
translations of coordinate axis. In this system of coordinates the above equations hold and $Z$ has the above form.
Moreover, the matrix $J$ becomes diagonal by blocks:
$$
J = \left(
\begin{array}{cll}
P & & 0\\
0 & & S
\end{array}
\right).
$$
where $P$ and $S$ are real $2 \times 2$-matrices
satisfying $P(0) = S(0) = J_{st}$.

In view of condition (\ref{tangency}) the Taylor expansion of $r$ is these coordinates has the form
\begin{eqnarray*}
\rho(w,z) = \Re z + \Re \sum_{j=1}^m q_j(w)z^j + p_m(w,\overline w) + o(\vert (w,z) \vert^m),
\end{eqnarray*}
where $p_m$ is a non-zero homogeneous polynomial of degree $m$ in $w$, $\overline w$ and
$q_1(0) = 0$.

\begin{e-lemme}
The polynomial $p_m$ contains at least one non-zero non-harmonic term.
\end{e-lemme}
\proof We need the following  statement due to Ivashkovitch-Rosay \cite{ivro}.

\begin{e-pro}
Let $m \geq 1$ integer  and $0 < \alpha < 1$. Let $J$ be a $C^{m+\alpha}$ almost
complex structure defined in a neighborhood of $0$ in $\R^{2n}$.
If $\phi: \D \longrightarrow (\R^{2n},J)$ is a smooth map such that
$\overline\partial_J \phi(\zeta) = o(\vert \zeta \vert^{m-1})$,
then there exists a $J$-holomorphic map $u$ from a neighborhood of $0$ in
$\cc$ into $\R^{2n}$ such that $\vert (u - \phi)(\zeta) \vert = o(\vert \zeta \vert^m)$.
\end{e-pro}

Suppose by contradiction that $p_m(w,\overline w) = ( h_m w^m + \overline{ h_m w^m})/2$.
Consider the map $\phi(\zeta) = (\zeta, -h_m \zeta^m)$.  In view of (\ref{vanishing})
it satisfies the hypothesis of the last proposition so there exists
a $J$-holomorphic disc $U $ of the form
$u = (\zeta, -h_m\zeta^m) + o(\vert \zeta \vert^m)$.
Then $\rho \circ u = o(\vert \zeta \vert^m)$
which contradicts the condition that the type of $E$ at the origin is equal to $m$.
This proves the lemma.
\bigskip

Now we are able to finish the proof. For $\delta > 0$ consider the non-isotropic dilations
$\Lambda_\delta: (w,z) \mapsto (\delta^{-1/m}w,\delta^{-1}z)$. Since in our system of coordinates $J$ is given
by a real $(4 \times 4)$-matrix diagonal by $(2 \times 2)$-blocks, the structures
$J_{\delta}: = (\Lambda_\delta)_*(J)$ tend to $J_{st}$ in  $C^\alpha$-norm for any positive $\alpha$.
The hypersurfaces $E_\delta : = \Lambda_\delta(E)$ tend to the  hypersurface
$\tilde E = \{ \Re z + p_m(w,\overline w) = 0 \}$ which is of finite type in view of the last lemma.
So we can apply Proposition \ref{perturbation}.

\bigskip

We also point out that a Bishop disc transversally attached to $\tilde E$ easily can be given quite explicitly
without using the general theorem of \cite{Tu}. Consider a holomorphic function of the form
$ w (\zeta) = e^{i\theta}(\zeta -1 )$ where a  value of the parameter $\theta \in [0,2\pi]$ will be chosen later.
After a standard biholomorphic change of coordinates we can assume that the polynomial $p_m$ contains
no harmonic terms, so that $p_m = \sum_{k=1}^{m-1}q_k w^k{\overline w}^{m-k}$ so that for
$\zeta \in \partial \D$ we have
$$
h(\zeta): = (p_m \circ w)(\zeta) = (\zeta -1)^{m} \sum_{k=1}^{m-1} (-1)^{m-k} e^{i(2k-m)\theta}q_k \zeta^{k-m}.
$$

Let $z$ be a holomorphic function on $\Delta$ such that
$\Re z \vert_{ \partial \D} = h$ and $z(1) = 0$.
The normal derivative $(\Re z)_{\nu}(1)$  at 1
(that is the derivative with respect to $\Re \zeta$) of $\Re z$
is given by the derivation of the Poisson integral (since $h(1) = 0$):
$$
(\Re z)_{\nu}(1) = (p.v.) \frac{1}{\pi}\int_0^{2\pi}
\frac{h(e^{i\tau})d\tau}{\vert e^{i\tau} - 1 \vert^2} =(p.v.)
\frac{i}{\pi} \int_{\zeta \in \partial \D}
\frac{h(\zeta)}{(\zeta - 1)^2}d\zeta.
$$
By the Cauchy residue theorem we obtain that
the last integral is equal to $\sum_{k=1}^{m-1} {q_k}\alpha_k e^{i(2k-m)\theta}$ where $\alpha_k$
are some non-zero constants coinciding up to the signs with certain binomial coefficients.
But this expression does not vanish after a suitable choice of $\theta$.
So for any $\varepsilon > 0$ the map $\zeta \mapsto (\varepsilon w(\zeta), \varepsilon^{m} z(\zeta))$
is an (arbitrarily small) Bishop disc transversally attached to $\tilde E$ at the origin.

If $E$ is a real analytic hypersurface in a real analytic almost complex manifold of complex dimension 2, then
 $E$ is of finite type if and only if it contains no
$J$-holomorphic discs ( by Nagano's theorem \cite{Na}).

This implies our main Theorem \ref{MainTheorem} in the real analytic category.

\begin{e-cor}
Let $(M,J)$ be a real analytic almost complex manifold of complex dimension 2 and $E$ be a real analytic
hypersurface in $M$. Assume that $E$ contains no $J$-holomorphic discs. Then the Bishop discs of $E$
fill a one-sided neighborhood of every point of $E$.
\end{e-cor}

A suitable modification of this method leads to a similar result in arbitrary dimension.
We say that a real hypersurface $E$ in an almost complex manifold $(M,J)$ is of finite type $m$
at point $p \in M$ if there exists a complex tangent vector $X$ to $E$ at $p$ such that any regular $J$-holomorphic
disc tangent to $X$ has the order of tangency with $E$ less or equal to $m$.

\begin{e-pro}
If $E$ is a real finite type hypersurface in an almost complex manifold $(M,J)$ (of an arbitrary dimension), then
its Bishop discs fill a one-sided neighborhood of every point of $E$.
\end{e-pro}

We can assume that local coordinates are chosen so that the disc
$Z = (\zeta,0,....,0)$ is $J$-holomorphic and
$J(0) = J_{st}$.

Then $ J = J_{st} + R$ where $R$ is $(2n \times 2n)$-matrix
formed by $2 \times 2$-blocks
$R_{kl}$ and $R(0) = 0$.
Moreover,
it follows
from the $J$-holomorphicity of $Z$ that
\begin{eqnarray}
\label{coefficients}
R_{k1}(\zeta,0,...,0) \equiv 0.
\end{eqnarray}
Furthermore, we can assume that $\rho(w^1,...,w^{n-1},z) = \Re z + o(\vert (w,z) \vert)$ and
$(\rho \circ Z)(\zeta) = p_m(\zeta,\overline \zeta) + o(\vert \zeta \vert^m)$.
Here $p_m$ is a homogeneous polynomial of degree $m$. Using the above Ivashkovich-Rosay proposition
and repeating the former argument we  obtain that $p_m$ contains a non-zero harmonic term,
because the condition (\ref{coefficients}) implies that $\overline\partial_J$ vanishes
with order $m-1$ on the disc $\zeta \mapsto (\zeta,0,....,0, a \zeta^m)$.

Finally, consider the dilations $\Lambda_\delta(z) = (\delta^{-1/m}w^1,\delta^{-1}w^2,...,\delta^{-1}z)$.
Again using condition (\ref{coefficients}),
we obtain that the structures $J_\delta = (\Lambda_\delta)_*(J)$ converge to the standard structure
$J_{st}$ in any $C^{\alpha}$ norm. The hypersurfaces $E_\delta = \Lambda_\delta(E)$ converge to the
hypersurface $E_0: \Re z + P_m(w^1,\overline w^1) = 0$ in $\cc^n$ which of finite
type with respect to $J_{st}$. So we conclude as above.

\section{Analysis of the Bishop equation}
Our goal now is to prove Theorem \ref{MainTheorem} stated in
introduction without any additional
hypothesis of pseudoconvexity or real analyticity type.
Everywhere we suppose that
$E = \{ \rho = 0 \}$ is a real $C^\infty$-smooth hypersurface
in an almost complex manifold $(M,J)$
of complex dimension 2. Since the statements of our main results
are local, we work in local coordinates similarly to the
previous sections. However, for technical reasons it is more
convenient to consider Bishop discs
attached to $E$ at the point $(0,1)$.
We also assume that $J(0,1) = J_{st}$.
Since in  the present section matrix computations will play
a substantial role, we everywhere write
vectors $Z$ of $\cc^2$ as vector-columns.
We will use the following notation:
$$
\overline\partial  := \frac{\partial}{\partial \overline \zeta},\;
\partial : = \frac{\partial }{\partial \zeta}.
$$
Recall that a map $\zeta \mapsto Z(\zeta)$ from
the unit disc $\D$ to $\cc^2$ is $J$-holomorphic if
and only if it satisfies the following system of equations:
\begin{eqnarray}
\label{CR'}
\overline\partial Z - A(Z) \overline{\partial Z} = 0,
\end{eqnarray}
where the matrix valued function $A$ is defined by (\ref{matrixA}).
We consider here only maps valued in a small enough
neighborhood of the point $(0,1)$.

Consider small embedded $J$-holomorphic Bishop discs
attached to $E$ at the point $(0,1)$;
their existence is proved in Section 3.
Given such a disc $Z_0$ there exists a
local diffeomorphism  such that in the new coordinates we have
$Z_0(\zeta) =  \left(
\begin{array}{cl}
0\\
\zeta
\end{array}
\right)$, $\zeta\in\D$
(see Lemma \ref{normalization}).
We denote again by $J$ the representation of our almost complex
structure in the new coordinates,
$J(0,1) = J_{st}$.

Our main goal is to establish the following

\begin{e-pro}
\label{MainProp}
Suppose that for the Bishop disc
$Z_0(\zeta) =  \left(
\begin{array}{cl}
0\\
\zeta
\end{array}
\right)$, $\zeta\in\D$,
we have $A \circ Z_0 = 0$ and $A_Z \circ Z_0 = 0$ and
\begin{eqnarray}
\label{norm}
\rho_z(0,1) = 1, \; \rho(0,1) = \rho_w(0,1) = 0.
\end{eqnarray}
Suppose that the derivatives $A_{\overline Z}$
and the second derivatives of $\rho$ are small
on $Z_0(\D)$.
Suppose further that for every Bishop disc
$Z(\zeta) =\left(
\begin{array}{cl}
z(\zeta)\\
w(\zeta)
\end{array}
\right)
$
with
$Z(1) = Z_0(1) =  \left(
\begin{array}{cl}
0\\
1
\end{array}
\right)
$
close enough to $Z_0$ we have $ \partial z(1) = 0$.
Then the Levi form of $E$ with respect to $J$
vanishes on the set $Z_0(b\D)$.
\end{e-pro}

Theorem \ref{keytheo} is an immediate corollary of this statement.
\bigskip

{\bf Remark 1.}
The assumptions of smallness of the derivatives of $A$ and
$\rho$ are automatically satisfied because the original
disc is small and the change of coordinates stretches it
to the unit disc.
The normalization condition (\ref{norm}) also can be achieved
by a suitable change of coordinates which does not affect previous
assumptions (Lemma \ref{normalization}).

{\bf Remark 2}. We point out that Proposition \ref{MainProp}
employs all $J$-holomorphic discs close to $Z_0$, not just
the infinitesimal perturbations of $Z_0$. In some sense it uses
the second variation of the condition of tangency.
Given a disc $Z_0$ of class $C^k$ for sufficiently big $k$,
we prove the desired vanishing of the Levi form under the hypothesis
that there are no $C^{k-m}$ transversal disc close to $Z_0$,
where $m$ is an unimportant constant. The ``loss of smoothness''
is due to the normalization of $A$ (Lemma \ref{normalization})
and changes of variables in the linearized Bishop equation
involving $A$ and $\rho$. It is unrelated to the class
$C^{1+\alpha}$, in which we consider the infinitesimal
perturbations in Section 5.3.

We provide a proof in the next three subsections.
In Section 5.1 we introduce and simplify the linearized Bishop
equation. In section 5.2 we describe its solutions for every
disc $Z$ close to $Z_0$. We use infinitesimal perturbations
$\dot Z$ of $Z$ to show that if no transversal perturbation exists,
then $Z$ satisfies a certain nonlinear integral equation $\Phi(Z)=0$
involving $Z$ and its derivatives. We only use the infinitesimal
perturbations $\dot Z$ parametrized by polynomials, hence they
are just as smooth as the coefficients of the linearized
Bishop equation, so no loss of smoothness occurs here.
Finally in Section 5.3 we use infinitesimal perturbations
of $Z_0$ of class $C^{1+\alpha}$ to prove the desired conclusion
about the Levi form.

\subsection{Linearized Bishop equation}
Suppose the hypotheses of Proposition \ref{MainProp}
are fulfilled. For simplicity we assume that
$A$ and $\rho$ are $C^\infty$.
Consider a Bishop disc
$Z(\zeta) =  \left(
\begin{array}{cl}
 z(\zeta)\\
 w(\zeta)
\end{array}
\right)
$, $\zeta\in\D$,
attached to $E$ at the origin and close enough to $Z_0$.
Then it satisfies the following boundary problem for the
$\overline\partial_J$-operator:
\begin{eqnarray*}
& &\overline\partial Z = A \overline{\partial Z},\\
& &\rho \circ Z \vert_{b\D} = 0.
\end{eqnarray*}
Since we only deal with small Bishop discs, we can assume
that for any fixed $k > 0$
the norm $\parallel A \parallel_{C^k}$ is small enough.

Consider now solutions $\dot Z$ of the corresponding linearized
problem:
\begin{eqnarray*}
& &\overline\partial\dot Z = \dot A \overline{\partial Z}
+ A \overline{\partial\dot Z},\\
& &\Re (\rho_Z(Z)\dot Z)\vert_{b\D} = 0.
\end{eqnarray*}
Here and in the rest of the paper for a map $F$ defined on the space
of smooth discs $Z :\D \longrightarrow \cc^2$, we denote by $\dot F$
the Frechet derivative of $F$.
In particular, if $F$ is defined by a smooth function on a set
in $\cc^2$, which is the case for $A$, then
\begin{eqnarray*}
\dot F(Z)(\dot Z) = F_Z(Z)\dot Z
+ F_{\overline Z}(Z)\overline{\dot Z}.
\end{eqnarray*}
We call the solutions $\dot Z$ of the linearized system
the {\it infinitesimal perturbations of the disc $Z$}.

In order to eliminate $\partial\dot Z$ in the right-hand side
of the linearized equation,
we perform the following change of (dependent) variables
$Z_1:= \dot Z - A \overline{\dot Z}$. Then
$\dot Z =  I' (Z_1 + A \overline Z_1)$,
where $ I' = (I - A \overline A)^{-1}$ and
$I$ denotes the identity matrix.

We have
\begin{eqnarray*}
& &\overline\partial Z_1=\overline\partial(\dot Z - A
\overline{\dot Z}) = (\overline\partial\dot Z-A
\overline\partial \overline{\dot Z})-\overline\partial A
\overline{\dot Z},\\
& &\overline\partial Z_1=\dot A
\overline\partial \overline Z - \overline\partial A
\overline{\dot Z},\\
& &\overline\partial A = A_Z \overline\partial Z
+A_{\overline Z} \overline\partial\overline{Z} = (A_ZA
+A_{\overline Z}) \overline\partial \overline{Z},\\
\end{eqnarray*}
Since $A_Z$ and $A_{\overline Z}$ are 3-index quantities,
the above products are not all matrix products, but rather
tensor products suitably contracted. We don't specify
the exact meaning of each product because it won't
matter for our analysis. We leave the details to the reader.
The same applies to several lines below. We get
\begin{eqnarray*}
& &\overline\partial Z_1 = A_1 Z_1 + A_2 \overline Z_1,\\
& &A_1 = (A_Z I' + A_{\overline Z} \overline{I'}\overline A)
\overline\partial \overline Z-(A_ZA + A_{\overline Z})
\overline\partial \overline Z \overline{I'}\overline A,\\
& &A_2 = (A_Z I' A + A_{\overline Z} \overline{I'})
\overline\partial \overline Z-(A_ZA + A_{\overline Z})
\overline\partial \overline Z \overline{I'}.
\end{eqnarray*}
Consider the Frechet derivatives $\dot A_1$ and $\dot A_2$
at $Z_0$. The condition $A = 0$ implies $\dot{I'} = 0$.
Obviously $\dot A_1$ and $\dot A_2$ are linear combinations
of $\dot Z$, $\overline{\dot Z}$
and $\overline\partial \overline{\dot Z}$ with smooth
coefficients depending on $\zeta$,
which we write in the form
\begin{eqnarray*}
\dot A_{\nu}(Z_0) = 0 \,\,\, \hbox{mod} \,
(\dot Z, \overline{\dot Z},
\overline\partial \overline{\dot Z}), \,\, \nu = 1,2.
\end{eqnarray*}
For the disc $Z=Z_0$ we have
$\overline\partial\overline Z=\overline\partial\overline Z_0=\left(
\begin{array}{cl}
0\\
1
\end{array}
\right).
$
The expression for $A_1$ implies
\begin{eqnarray*}
\dot A_1(Z_0) = 0 \,\,\, \hbox{mod} \, ( \dot Z, \overline{\dot Z}).
\end{eqnarray*}

We now express the boundary condition $\Re ( \rho_Z \dot Z )=0$
in terms of $Z_1$. We obtain
\begin{eqnarray*}
\rho_Z I' (Z_1 + A \overline Z_1) +
\rho_{\overline Z} \overline{I'}(\overline Z_1 + \overline A Z_1)=0,
\end{eqnarray*}
We put
\begin{eqnarray*}
\lambda =(\lambda_1,\lambda_2) = \rho_Z I'
+ \rho_{\overline Z} \overline{I'}\overline A.
\end{eqnarray*}
Then the boundary condition for $Z_1$ turns into
\begin{eqnarray*}
\Re (\lambda Z_1)|_{b\D}=0.
\end{eqnarray*}

For simplicity of notation we put
$$
\partial_1 f = \partial f(1).
$$

We note that the condition $\partial_1 z = 0$
fulfilled for all the discs implies
$\partial_1 \dot z = 0$.
For the change of variables
$\dot Z \mapsto Z_1 =   \left(
\begin{array}{cl}
z_1\\
w_1
\end{array}
\right)$, the condition $A(0,1) = 0$ (recall $J(0,1) = J_{st}$ )
implies that
$\partial_1 z_1 = \partial_1 \dot z$.
Thus for every $Z_1$ the conditions
\begin{eqnarray*}
& &\overline\partial Z_1 = A_1 Z_1 + A_2 \overline Z_1,\\
& &Z_1(1) = 0,\\
& &\Re(\lambda Z_1)|_{b\D} = 0
\end{eqnarray*}
imply
$\partial_1 z_1 = 0$.
Consider the matrix
$$\Lambda =  \left(
\begin{array}{cll}
\lambda_1 & & \lambda_2\\
0 & & 1
\end{array}
\right),
$$
where $\lambda = (\lambda_1, \lambda_2)$ is defined above.
We use that $\Lambda$ is smoothly extended
on the disc $\D$, although only the values on $b\D$ will matter.

In order to simplify the boundary conditions, we further
change the variable to $V = \Lambda Z_1$ with
$V =  \left(
\begin{array}{cl}
v_1\\
v_2
\end{array}
\right)$.
We put
\begin{eqnarray*}
& &B_1 = (\overline\partial \Lambda)\Lambda^{-1}
+ \Lambda A_1 \Lambda^{-1},\\
& &B_2 = \Lambda A_2 \overline{\Lambda^{-1}}.
\end{eqnarray*}
It follows from the hypotheses of Proposition \ref{MainProp}
that $B_1$ and $B_2$ are small on $Z_0$, which we will
need in Section 5.2.
The new unknown $V$ satisfies
\begin{eqnarray*}
& &\overline\partial V = B_1 V + B_2\overline V,\\
& &V(1)=0,\\
& &\Re v_1|_{b\D} =  0.
\end{eqnarray*}
For every solution $V$ we have
$$\partial_1 v_1 = 0.$$

The expressions for $\dot A_1$ and $\dot A_2$ on the disc $Z_0$ imply
\begin{eqnarray}
\label{B1dot}
& &\dot B_1 = (\overline\partial \dot \Lambda) \Lambda^{-1} \,\,\,
\hbox{mod} \, (\dot Z, \overline{\dot Z}),\\
\label{B2dot}
& &\dot B_2 = 0 \,\,\,
\hbox{ mod} \, (\dot Z, \overline{\dot Z}, \overline\partial
\overline{\dot Z}).
\end{eqnarray}
We will need this information about the derivatives
$\dot B_\nu$ in Section 5.3.

\subsection{Solving the boundary value problem: the generalized
Schwarz integral}
We derive an integral
formula for solving a linear boundary Riemann-Hilbert problem
for elliptic systems of PDE in the unit disc.
Following \cite{Ve}, we call it the
generalized Schwarz integral.
Consider the following boundary value problem in $\D$:
\begin{eqnarray*}
& &\overline\partial V = B_1 V + B_2 \overline V,\\
& &\Re V|_{b\D} = \Gamma \,\,\,(\Gamma(1) = 0),\\
& &V(1) = 0.
\end{eqnarray*}
Recall the    Cauchy integral $K$ and the Cauchy-Green
transform $T$
are defined in the previous sections by (\ref{Cauchy})
and (\ref{CauchyGreen}).
We also define
the operators:
\begin{eqnarray*}
K_1u: = Ku - (Ku)(1)
\end{eqnarray*}
and
\begin{eqnarray*}
T_1u: = Tu - (Tu)(1).
\end{eqnarray*}
We also need the following operators:
\begin{eqnarray*}
T_1^*u = T^*u - (T^*u)(1),
\end{eqnarray*}
where
\begin{eqnarray*}
T^*u(\zeta) = \frac{1}{2 \pi i} \int\int_{\D}
\frac{\zeta u(\tau)d\tau \wedge d\overline\tau}{1-\zeta\bar\tau}.
\end{eqnarray*}
The operator $T^*$ satisfies the following identities:
\begin{eqnarray*}
T^* = -K \overline{T}, \;
T_1^* = -K_1\overline{T}_1,
\end{eqnarray*}
which one can verify directly.

Then by the  Schwarz-Green formula:
\begin{eqnarray*}
V = 2K_1 \Re V + T_1 \overline\partial V
+ T_1^*\overline{\overline\partial V} + V(1).
\end{eqnarray*}
Set $V_0 = 2 K_1 \Gamma$. Then $V$ satisfies the following
integral equation
\begin{eqnarray*}
V = V_0 + T_1(B_1 V + B_2 \overline V) +
T_1^*(\overline B_1 \overline V + \overline B_2 V).
\end{eqnarray*}
Consider the operators
\begin{eqnarray*}
& &L_1 = T_1 \circ B_1 + T_1^* \circ \overline B_2,\\
& &L_2 = T_1 \circ B_2 + T_1^* \circ \overline B_1,
\end{eqnarray*}
where  $B_j$ are viewed as operators of left multiplication
by the matrix functions $B_j$.
Then $V$ satisfies the following equation:
\begin{eqnarray*}
V = V_0 + L_1 V + L_2 \overline V.
\end{eqnarray*}
We solve this equation using that $L_1$ and $L_2$ are small
because so are $B_1$ and $B_2$.
Then the solution has the form
\begin{eqnarray*}
V = V_0 + R_1V_0 + R_2\overline V_0,
\end{eqnarray*}
where the resolvent operators $R_j$ are $\cc$-linear.
We now find an explicit expression for the resolvents $R_j$.
We have
\begin{eqnarray*}
V_0 + R_1 V_0 + R_2 \overline V_0 = V_0
+ L_1(V_0 + R_1 V_0 + R_2 \overline V_0) +
L_2(\overline V_0 + \overline R_1 \overline V_0
+ \overline R_2 V_0).
\end{eqnarray*}
Then (we temporarily drop the index $0$ of $V_0$)
\begin{eqnarray*}
& &R_1V = L_1(V + R_1V) + L_2\overline R_2 V,\\
& &R_2\overline V = L_1 R_2 \overline V + L_2
(\overline V + \overline R_1 \overline V),
\end{eqnarray*}
so that
\begin{eqnarray*}
& &R_1 = L_1 + L_1R_1 + L_2\overline R_2,\\
& &R_2 = L_1 R_2 + L_2(I + \overline R_1).
\end{eqnarray*}
Substituting the equality
\begin{eqnarray*}
R_2 = (I - L_1)^{-1}L_2(I + \overline R_1)
\end{eqnarray*}
into the expression for $R_1$, we obtain
\begin{eqnarray*}
R_1 = [I - L_1 - L_2(I- \overline L_1)^{-1}\overline L_2]^{-1}
[L_2(I - \overline L_1)^{-1}\overline L_2 + L_1].
\end{eqnarray*}
Since
$(1-x)^{-1}x=(1-x)^{-1}-1$,
we can rewrite the above expression in the form
\begin{eqnarray*}
R_1 = [I - L_1 - L_2(I - \overline L_1)^{-1}\overline L_2]^{-1} - I.
\end{eqnarray*}
Hence $R_j$ are bounded operators $C^m \longrightarrow C^{m+1}$
for any non-integral $m > 0$.
We use the expression for $R_1$ to make some simple
observations. Note that
\begin{eqnarray*}
R_1 = L_1 + O(2),
\end{eqnarray*}
where $O(2)$ is a sum of products of 2 or more operators
$L_1$, $L_2$, $\overline L_1$, $\overline L_2$.
It will turn out that the $O(2)$ term is unimportant,
but it will take considerable work. Furthermore, one can see that
every term in $R_1$ starts with $L_1$ or $L_2$ and ends
with $L_1$ or $\overline L_2$.
Recalling the expressions for $L_1$ and $L_2$
we conclude that $R_1$ is a sum of terms of the form
\begin{eqnarray*}
P_0 \circ F_0 \circ P_1 \circ F_1 \circ ....\circ P_n \circ F_n,
\end{eqnarray*}
where each
$P_j\in\{T_1, \overline T_1, T_1^*, \overline T_1^*\}$ and
$F_j\in\{B_1, \overline B_1, B_2, \overline B_2\}$.
Moreover,
$P_0\in\{T_1, T_1^*\}$ and
$F_n\in\{B_1, \overline B_2\}$.

We now interpret the condition $\partial_1 v_1:= \partial v_1(1)= 0$.
The boundary conditions we derived in Section 5.1 imply
$\Gamma =  \left(
\begin{array}{cl}
0\\
\gamma
\end{array}
\right)$,
where $\gamma$ is an arbitrary real function with
$\gamma(1)=0$. Then
$V_0 =  \left(
\begin{array}{cl}
0\\
\varphi
\end{array}
\right)$,
where $\varphi$ is an arbitrary holomorphic function
$\varphi = 2K_1 \gamma$.
Then $v_1=R_1^{12}\varphi+R_1^{12}\overline\varphi$,
where $R_1^{12}$ and $R_2^{12}$ are the (12) matrix
entries of the matrix operators $R_1$ and $R_2$.
The condition
$\partial_1 v_1 = 0$ means
\begin{eqnarray*}
\partial_1 R_1^{12} \varphi + \partial_1 R_2^{12}
\overline \varphi=0.
\end{eqnarray*}
Since the first term is $\cc$-linear in $\varphi$ and
the second one is antilinear, we get
\begin{eqnarray*}
\partial_1 R_1^{12}\varphi = 0
\end{eqnarray*}
for all $\varphi$ with $\varphi(1) = 0$.

Since every term in $R_1$ starts with $T_1$ or $T_1^*$,
we need formulas for their derivatives at $\zeta=1$.
For every function $f$ of class $C^\alpha(\overline\D)$,
$0<\alpha<1$, with $f(1) = 0$ we have
\begin{eqnarray*}
& &\partial_1 T_1f = \partial_1 Tf= \frac{1}{2\pi i}
\int\int_{\D} \frac{f(\tau)
d\tau \wedge d\overline\tau}{(\tau - 1)^2},\\
& &\partial_1 T_1^*f = \partial_1 T^*f
= \frac{1}{2 \pi i}\int\int_{\D} \frac{f(\tau)
d\tau \wedge d\overline \tau}{(\overline\tau - 1)^2},
\end{eqnarray*}
where the integrals converge in the usual sense.

We would like to interpret the condition
$\partial_1 R_1^{12}\varphi = 0$ for all $\varphi$
as a moment condition and eliminate the arbitrary
$\varphi$. In order to do so we need to reverse the
order of integration in every term in $R_1$, which
requires the following preparations.
Let $P$ be a scalar integral operator with
kernel $P(t,\tau)$, that is
\begin{eqnarray*}
Pf(\zeta) = \int\int_{\D} P(\zeta,\tau) f(\tau)
d\tau \wedge d\overline \tau.
\end{eqnarray*}
Define the operators $P^+$ and $P^-$ as the integral operators
with the kernels
\begin{eqnarray*}
& &P^+(\zeta,\tau) = \frac{(\zeta-1)^2}
{(\tau - 1)^2}P(\tau,\zeta),\\
& &P^-(\zeta,\tau) = \frac{(\overline \zeta -1)^2}
{(\overline \tau - 1)^2}P(\tau,\zeta).
\end{eqnarray*}
We also put
$$\mu(\tau) = \left ( \frac{\tau - 1}
{\overline \tau - 1} \right )^2$$

\begin{e-lemme}
\label{Lemma1}
Let $\varphi(\tau) = (\tau - 1)^2 \psi(\tau)$.
We have
\begin{eqnarray*}
\partial_1TF_0 \circ P_1 \circ F_1 \circ ...\circ P_n
\circ F_n\varphi=\frac{1}{2 \pi i}\int\int_{\D}
P^+_n ( ... P^+_1(F_0)F_1 ... ) F_n \psi(\tau)
d\tau \wedge d\overline\tau,
\end{eqnarray*}
\begin{eqnarray*}
\partial_1T^*F_0 \circ P_1 \circ F_1 \circ ... \circ P_n
\circ F_n\varphi=\frac{1}{2\pi i} \int\int_{\D}
P^-_n ( ... P^-_1(F_0)F_1 ... ) F_n \psi(\tau) \mu(\tau)
d\tau \wedge d\overline\tau.
\end{eqnarray*}
\end{e-lemme}
The proof follows by changing the
order of integration.

\begin{e-lemme}
\label{Lemma2}
(Moment conditions.)
Let $F \in L^\infty(\D)$. Then
$$
TF\vert_{b\D} = 0
$$
if and only if for every holomorphic polynomial $\varphi$ we have
$$
\int\int_{\D}F(\tau)\varphi(\tau)d\tau \wedge d\overline\tau = 0.
$$
\end{e-lemme}
\proof
The series
$$
\frac{1}{\tau-\zeta}=-\sum_{n=0}^\infty\frac{\tau^n}{\zeta^{n+1}}
$$
converges in $L^1(\D)$ for every fixed $\zeta\in b\D$.
Integrating the series against $F(\tau)$ yields
\begin{eqnarray*}
TF(\zeta) = -\sum_{n=0}^{\infty} \frac{c_n}{\zeta^{n+1}},\;
\zeta\in b\D,\;\;\;
\hbox{where}\;\;\;
c_n = \frac{1}{2\pi i}\int\int_{\D}F(\tau)\tau^n
d\tau \wedge d\overline\tau,
\end{eqnarray*}
and the lemma follows.
\bigskip

\begin{e-lemme}
\label{Lemma5}
Let $f \in C^\gamma(\overline \D)$, $0<\gamma < 1$,
$f(1) = 0$. Let $\vert g \vert \leq C$ be bounded and
\begin{eqnarray*}
\vert \partial g(\tau) \vert \leq
\frac{C}{\vert \tau - 1 \vert},\;\;
\vert \overline\partial g(\tau) \vert \leq
\frac{C}{\vert \tau - 1 \vert},\;\;
\tau \in \overline \D \backslash \{ 1 \}.
\end{eqnarray*}
Then $f g \in C^\gamma(\overline \D)$ and
$\parallel f g \parallel_{C^\gamma} \leq
C'C\parallel f \parallel_{C^\gamma}$ for some constant $C'$.
For $\gamma = 1$ lemma holds in $\Lip_1(\overline\D)$.
\end{e-lemme}
\proof
Put
$\Delta:= (gf)(\zeta_1) - (gf)(\zeta_2) = (f(\zeta_1) - f(\zeta_2))
g(\zeta_1) + f(\zeta_2)(g(\zeta_1) - g(\zeta_2))$,
$\vert \zeta_1 - \zeta_2 \vert = \delta$.
Without loss of generality we can assume that
$\vert \zeta_2 -1 \vert \leq \vert \zeta_1 - 1 \vert$.
Then
\begin{eqnarray*}
\vert \Delta \vert \leq
C \parallel f \parallel_{C^\gamma}\delta^\gamma
+\parallel f \parallel_{C^\gamma}\vert \zeta_2
- 1 \vert^\gamma \min \left ( 2C,\frac{\pi C \delta}
{\vert \zeta_2 - 1 \vert} \right ) \leq \pi \delta,
\end{eqnarray*}
where $\pi$ arises because we may have to connect
$\zeta_1$ and $\zeta_2$ by going around 1.

Now consider the two cases separately:
$\delta \leq \vert \zeta_2 - 1 \vert$ and $\delta \geq
\vert \zeta_2 - 1 \vert$.
In both cases we get
$\vert \Delta \vert \leq C'C\parallel f \parallel_{C^\gamma}
\delta^\gamma$,
which proves the lemma.
\bigskip

\begin{e-lemme}
\label{Lemma3}
If $P$ is one of the operators $T_1$, $T_1^*$,
$\overline T_1$, $\overline T_1^*$, then $P^+$ and $P^-$
are bounded operators
$C^\alpha(\overline\D) \longrightarrow C^{1-\beta}(\overline\D)$
for every $0 < \alpha < 1$, $0 < \beta < 1$.
\end{e-lemme}
\proof We directly obtain the following equalities:
\begin{eqnarray}
\label{1'}
T_1^+ = -T_1,
\end{eqnarray}
\begin{eqnarray}
\label{2'}
T_1^- = -\mu^{-1}T_1\mu,
\end{eqnarray}
\begin{eqnarray}
\label{3'}
{T_1^*}^+ = -\mu\overline{T_1^*}.
\end{eqnarray}

By using the equality $T^* = -K\overline{T}$, we obtain
\begin{eqnarray}
\label{5'}
{T_1^*}^+ = \mu\overline K_1 T_1,
\end{eqnarray}
\begin{eqnarray}
\label{6'}
{T_1^*}^- = \overline K_1 T_1 \mu.
\end{eqnarray}

We also need the following relation:
\begin{eqnarray}
\label{7'}
T_1\mu(\zeta) = T\mu(\zeta) = \frac{\zeta -1}{\overline \zeta - 1}(1 - \vert \zeta \vert^2)
\in \Lip_1(\overline\D).
\end{eqnarray}
To obtain (\ref{7'}) we put
$\psi(\zeta) = -\frac{(\zeta - 1)^2}{\overline{\zeta} -1}$,
then $\overline\partial \psi = \mu$.
By the Cauchy-Green formula we have
$\psi = K\psi + T\mu$. Note that
$\psi(\zeta)\vert_{b\D} = \zeta(\zeta -1)$
so that $K\psi = \zeta(\zeta-1)$. Hence
$$
T\mu = \psi - K\psi = \frac{\zeta -1}{\overline \zeta  -1}
(1 - \vert \zeta \vert^2).
$$
Since $(\overline P)^+ = \overline{P^-}$,
$(\overline P)^- = \overline{P^+}$, it suffices to prove
the lemma for $P \in \{ T_1, T_1^* \}$.

For $T_1^+$ by (\ref{1'}) we have
$T_1^+:C^\alpha(\overline \D) \longrightarrow
C^{1+\alpha}(\overline \D)$.
Consider $T_1^-$. Let $u \in C^\alpha(\overline\D)$.
Then $T_1(\mu u) = T_1(\mu u(1)) + T_1(\mu(u - u(1)))$.
Then $T_1(\mu u(1)) = u(1) T_1\mu \in
\Lip_1(\overline \D)$ by (\ref{7'}).
Then $\mu(u - u(1)) \in C^\alpha(\overline\D)$
by Lemma \ref{Lemma5} and
$T_1(\mu(u - u(1)))  \in C^{1+\alpha}(\overline\D)$.
Therefore $T_1(\mu u) \in \Lip_1(\overline \D)$ and
$T_1^- u = -\mu^{-1} T_1(\mu u) \in \Lip_1(\overline\D)
\subset C^{1-\beta}(\overline\D)$
by Lemma \ref{Lemma5}.
The operators ${T_1^*}^+$ and ${T_1^*}^-$ are treated
similarly by (\ref{5'}) and (\ref{6'}).
This completes the proof of the lemma.
\bigskip

The property $\partial_1 R^{12}\varphi = 0$ for all holomorphic
functions $\varphi$ with $\varphi(1) = 0$, takes the form
\begin{eqnarray*}
\partial_1 [ T_1\circ B_1 + T_1^*\circ \overline B_2 +
\sum_{n \geq 1} P_0\circ F_0\circ ... \circ P_n\circ F_n ]^{12}
\varphi = 0,
\end{eqnarray*}
where the summation includes all the terms in $R_1$, so
$P_0,...,P_n \in \{ T_1, T_1^*, \overline T_1, \overline T_1^* \}$,
$P_0 \in \{ T_1, T_1^* \}$,
$F_n \in \{ B_1, \overline B_2 \}$ and $(12)$ denotes the
corresponding matrix entry.
By Lemma \ref{Lemma1} this condition is equivalent to
\begin{eqnarray*}
\int\int_{\D} \left [ B_1 + \overline B_2 \mu(\tau) + \Sigma' +
\Sigma'' \right ]^{12} \psi(\tau) d\tau \wedge d\overline\tau = 0,
\end{eqnarray*}
where
\begin{eqnarray*}
\Sigma' =\sum_{n \geq 1} P_n^+( ... P_1^+(F_0)F_1 ... )F_n, \qquad
\Sigma''=\sum_{n \geq 1} P_n^-( ... P_1^-(F_0)F_1 ... )F_n \mu(\tau),
\end{eqnarray*}
and $\psi$ is an arbitrary holomorphic function.
The terms in $\Sigma'$ and $\Sigma''$ correspond to $O(2)$
in $R_1$. They split into $\Sigma'$ and $\Sigma''$
depending whether the term in $R_1$ starts with
$T_1$ or $T_1^*$.
The series converges in $C^{\1-\beta}(\overline\D)$.

By Lemma \ref{Lemma2} we further obtain
\begin{eqnarray}
\label{functional}
\Phi(Z):= T\left [ B_1 + \overline B_2 \mu(\tau) + \Sigma'
+\Sigma'' \right ]^{12} \vert_{b\D} = 0.
\end{eqnarray}
This property holds for every Bishop disc $Z$ close enough
to $Z_0$ with
$Z(1) =  \left(
\begin{array}{cl}
0\\
1
\end{array}
\right)$
provided that for every polynomial $\psi$ the corresponding
infinitesimal perturbation $\dot Z$ of $Z$ is tangent to $E$
at $\zeta=1$.

\subsection{Frechet derivative of  $\Phi$}
In this section we consider the map
$\Phi$ defined by  (\ref{functional}).
This map is defined in a neighborhood of the disc $Z_0$
in the Banach manifold ${\cal A}_p^J(E)$ of $J$-holomorphic Bishop
discs attached to $E$ at $p=(0,1)$.
The condition (\ref{functional}) means
that in fact the map $\Phi$ vanishes identically,
so its Frechet derivative $\dot \Phi$ at $Z_0$ does.
We will study the  geometric consequences of the equality
\begin{eqnarray*}
\dot \Phi(Z_0)(\dot Z) = 0,\,\,\,\dot Z\in T_{Z_0}{\cal A}_p^J(E).
\end{eqnarray*}
According to the previous subsection, due to the special
normalization of the matrix $A$ along the disc $Z_0$, we have
$$
\dot Z = \Lambda^{-1}V
$$
$$
V = V_0 + R_1V_0 + R_2\overline V_0,\,\,\,\,
V_0 =  \left(
\begin{array}{cl}
0\\
\varphi
\end{array}
\right),
$$
where $\varphi$ is an arbitrary holomorphic function in $\D$
with $\varphi(1)=0.$

This allows to consider $\dot Z$ and therefore $\dot \Phi$
as $\R$-linear operators applied to
a function
$\varphi \in {\cal O}(\D) \cap C^{1+\alpha}
(\overline \D)$, $0 < \alpha < 1$.

For the target spaces, given integral $m \geq 0$ and
$0 < \alpha < 1$ we introduce the spaces
$C^{m + \alpha}_1(\overline\D)$ and $C^{m+\alpha}_1(b\D)$
as spaces of functions which are
$C^{m + \alpha}$ except at $1$.
More precisely we say that
$f \in C_1^{m+\alpha}(\overline\D)$ if $f \in L^\infty(\D)$
and for every $\varepsilon > 0$ we have
$f\vert_{\overline\D \backslash \B(1,\varepsilon)}
\in C^{m +\alpha}(\overline\D \backslash \B(1,\varepsilon))$,
where $\B(1,\varepsilon)$ denotes the disc of radius
$\varepsilon$ centered at $1$.
We define $C_1^{m+\alpha}(b\D)$ similarly.
We do not introduce norms in these spaces, nevertheless,
we say that $P$ is {\it a bounded linear operator}
$C^{m+\alpha}_1(\overline\D) \longrightarrow
C_1^{k+\beta}(\overline\D)$
if for every $\varepsilon > 0$ there exists a constant
$C = C(\varepsilon)>0$ such that for every
$f \in C^{m+\alpha}_1$ we have
\begin{eqnarray*}
\parallel Pf \parallel_{C^{k+\beta}(\D \backslash
\B(1,2\varepsilon))}
+\parallel P f \parallel_{L^\infty(\D)} \leq
C (\parallel f \parallel_{C^{m+\alpha}(\D \backslash
\B(1,\varepsilon))}
+\parallel f \parallel_{L^\infty(\D)}).
\end{eqnarray*}
Similarly, we define bounded operators
$C^{m+\alpha}(\overline\D) \longrightarrow
C_1^{k+\beta}(\overline\D)$ and
when we have $b\D$ in place of $\D$.

\begin{e-lemme}
\label{Lemma3'}
If $P \in \{ T_1, T_1^*, \overline T_1, \overline T_1^* \}$
then $P^+$ and $P^-$ are bounded operators
$C^{m+\alpha}_1(\overline \D) \longrightarrow
C^{m+1+\alpha}_1(\overline\D)$ for all integral $m \geq 0$,
$0 < \alpha <1$.
\end{e-lemme}
The lemma follows immediately by splitting integration
\begin{eqnarray*}
\int\int_{\D} = \int\int_{\D \backslash \B(1,\varepsilon)}
+ \int\int_{\B(1,\varepsilon)}
\end{eqnarray*}
in the definition of each operator and the regularity of
$T$ and $T^*$.
\bigskip

The Frechet derivative of $\Phi$ has the form
$$
\dot\Phi = T\left [ \dot B_1 + \overline{\dot B}_2 \mu(\tau)
+ \dot\Sigma' + \dot\Sigma'' \right ]^{12} \vert_{b\D}=0.
$$
In the derivatives $\dot\Sigma'$ and $\dot\Sigma''$
we distinguish the expressions
\begin{eqnarray*}
\Psi' = \sum_{n \geq 1} P_n^+ ( ... P_1^+ (F_0) ... )\dot{} \, F_n
\qquad\hbox{and}\qquad
\Psi''= \sum_{n \geq 1} P_n^- ( ... P_1^- (F_0) ... )\dot{} \, F_n,
\end{eqnarray*}
in which $F_n$ is not differentiated.
Since $F_n \in \{ B_1, \overline B_2 \}$,
then all other terms contain either
$\dot B_1$ or $\overline{\dot B}_2$, and we group them
accordingly. Then
$$
\dot \Phi = T \left [ (I + a_1 + a_2\mu)\dot B_1
+ (b_1 + b_2\mu)\overline{\dot B}_2 + \Psi' +
\mu \Psi'' \right ]^{12}\vert_{b\D} = 0.
$$
Here each matrix function $a \in \{ a_1, a_2, b_1, b_2 \}$
has the form
\begin{eqnarray*}
a = \sum P_n^\sigma
(P_{n-1}^\sigma(... P_1^\sigma (F_0) ...)F_{n-1}),
\end{eqnarray*}
where $\sigma$ is $+$ or $-$.
(In each term in $a$, the signs $\sigma$ are the same,
but different terms have different $\sigma$, depending
whether they come from $\Sigma'$ or $\Sigma''$.)
Then
$a\in C^{1-\beta}$ and the norms of $a_j$, $j=1,2$, are small.

Furthermore,  actually
$a\in C_1^{m+\alpha}(\overline\D)$ for all $m \geq 0$,
$0 < \alpha < 1$.
Indeed, since $F_j$ are smooth, then each term in $a$
is $C_1^{m+\alpha}$
for all $m$ and $\alpha$ by Lemma \ref{Lemma3'}.
To show that $a \in C_1^{m+\alpha}$ we split the terms
containing $m$ or more operators $P_j^\sigma$
into finitely many groups of the form
\begin{eqnarray*}
P_n^\sigma ( P_{n-1}^\sigma ( ... P_{n-m+1}^\sigma (\sum ...)
F_{n-m+1} ...) F_{n-1}).
\end{eqnarray*}
Each group consists of the terms with the same
$P_n^\sigma, ..., P_{n-m+1}^\sigma$ and $F_{n-m+1}, ..., F_{n-1}$.
Then the inside summation is of class
$C^{1-\beta}$ by Lemma \ref{Lemma3}, hence each
group is $C_1^{m + 1 -\beta}$ by Lemma \ref{Lemma3'}.
We say that $a \in C_1^\infty(\overline\D)$.

Since $F_j \in \{ B_1,B_2,\overline B_1, \overline B_2 \}$ and
$\dot Z \in C^{1+\alpha}$,
we have $\dot F_j \in C^\alpha$ (because
the expressions for $\dot B_1$, $\dot B_2$ include the term
$\overline\partial \overline{\dot Z}$).

By Lemma \ref{Lemma3'} the terms $\Psi'$ and
$\Psi''$ define bounded operators
$C^{1+\alpha}(\overline\D) \longrightarrow
C^{1+\alpha}_1(\overline\D)$.
Hence $T\Psi'$ defines a bounded operator
$C^{1+\alpha}(\overline\D) \longrightarrow
C^{2+\alpha}_1(\overline\D)$.
The same is true for $T(\mu\Psi'')$ because
the multiplication by $\mu$ does not change
the class $C^{m+\alpha}_1(\overline\D)$.

The next step is the following

\begin{e-lemme}
\label{B2bar}
The term $T[(b_1 + b_2\mu)\overline{\dot B_2}]$ defines
a bounded linear operator
$C^{1+\alpha}(\overline\D) \longrightarrow C^{2 + \alpha}_1(b\D)$.
\end{e-lemme}
This simple result applies to all matrix elements, including (12)-entry.
By (\ref{B2dot})
\begin{eqnarray*}
\overline{\dot B_2} = 0 \,\,\,
\hbox{mod} \,(\dot Z, \overline{\dot Z}, \partial \dot Z).
\end{eqnarray*}
The operator $T$ takes the terms with $\dot Z$ and
$\overline{\dot Z}$ to $C^{2+\alpha}_1(\overline\D)$.
Now since $\dot Z = \Lambda^{-1} V_0 + (S)$, where $(S)$
denotes smoother terms,
it suffices to consider $T(b \partial \varphi)$,
where $b \in C^\infty_1$.
We have
\begin{eqnarray*}
T(b\partial\varphi) = [T,b]\partial\varphi + bT(\partial\varphi),
\end{eqnarray*}
where $[T,b]:= T \circ b - b \circ T$ denotes the commutator
of two operators. In order to complete the proof of Lemma \ref{B2bar},
we need the following two simple lemmas.
\begin{e-lemme}e
\label{Lemma6}
Let $b\in C^\infty_1(\overline\D)$.
Then $[T,b]$ defines a bounded operator
$C^m(\overline \D) \longrightarrow C^{m+2}_1(\overline \D)$
for any non-integral $m > 0$.
\end{e-lemme}
The proof is immediate because the kernel of the operator $[T,b]$
has the form
\begin{eqnarray*}
\frac{b(\tau) - b(t)}{\tau - t}.
\end{eqnarray*}
\begin{e-lemme}
\label{Lemma7}
Let $\varphi$ be a holomorphic function. Then
\begin{eqnarray*}
T\varphi(\tau) = \overline \tau \varphi(\tau)
-\frac{\varphi(\tau) - \varphi(0)}{\tau}.
\end{eqnarray*}
In particular if $\vert \tau \vert = 1$,
then $T\varphi(\tau) = \frac{\varphi(0)}{\tau}$.
\end{e-lemme}
\proof Set $f = \overline\tau \varphi$. By the Cauchy-Green formula
$$
f = Kf + T\overline\partial f = Kf + T\varphi.
$$
Then $Kf = K(\varphi/\tau) = \frac{\varphi - \varphi(0)}{\tau}$
which proves the lemma.
\bigskip

We now conclude the proof of Lemma \ref{B2bar}.
It follows by Lemma \ref{Lemma6} and \ref{Lemma7}
that $[T,b]\partial \varphi \in C^{2 +\alpha}_1(\D)$ and
$T(\partial \varphi) \in C^{\infty}(b\D)$.
This proves Lemma \ref{B2bar}.
\bigskip

We will write $P_1 \sim P_2$ for two operators
$P_1$ and $P_2$ if $P_1 - P_2$ is a bounded operator
$C^{1+\alpha}(\overline\D) \cap {\cal O}(\D) \longrightarrow
C_1^{2+\alpha}(b\D)$.
The result of the above analysis of $\dot\Phi$
so far yields
$$
T((I + a_1 + a_2 \mu)\dot B_1)^{12}\vert_{b\D}\sim 0.
$$
Put $I_0 = I + a_1 + a_2 \mu$. Using the commutator argument
and Lemma \ref{Lemma6}, we obtain
$(I_0 T \dot B_1)^{12}\vert_{b\D}\sim 0$.
Then
$$
0\sim(I_0 T \dot B_1)^{12}\vert_{b\D}
=I_0^{11} T\dot B_1^{12}+I_0^{12} T\dot B_1^{22}.
$$
Since
$$
\dot \Lambda =  \left(
\begin{array}{cll}
\dot \lambda_1 & &\dot \lambda_2\\
0 & & 0
\end{array}
\right),
$$
then $T\dot B_1^{22}\vert_{b\D}\sim 0$.
Hence
$I_0^{11}(T\dot B_1)^{12}\vert_{b\D}\sim 0$.
Since $I_0$ is close to $I$, then
$$
(T\dot B_1)^{12}\vert_{b\D}\sim 0.
$$
Recall by (\ref{B1dot})
\begin{eqnarray*}
\dot B_1 = (\overline\partial \dot \Lambda)\Lambda^{-1} \,\,\,
\hbox{mod} \,
(\dot Z, \overline{\dot Z}).
\end{eqnarray*}
On the disc $Z_0$ we have
$$
\Lambda =  \left(
\begin{array}{cll}
\rho_z & &\rho_w\\
0 & &1
\end{array}
\right),
\qquad
\Lambda^{-1} =  \rho_z^{-1}
\left(
\begin{array}{cll}
1 & & -\rho_w\\
0 & &\rho_z
\end{array}
\right),
$$
$$
[(\overline \partial \dot
\Lambda)\Lambda^{-1}]^{12}=\rho_z^{-1}\overline\partial\dot\lambda
\left(
\begin{array}{cl}
-\rho_w\\
\rho_z
\end{array}
\right).
$$
By Lemma \ref{Lemma6},
$T(\dot B_1)^{12} \vert_{b\D} \sim 0$
implies
\begin{eqnarray}
\label{Tdbar}
T(\overline\partial\dot\lambda)\vert_{b\D}
\left(
\begin{array}{cl}
-\rho_w\\
\rho_z
\end{array}
\right)
\sim 0
\end{eqnarray}
Since $A_Z\circ Z_0=0$, then
$$
\dot \lambda = (\rho_Z)\dot{} + \rho_{\overline Z}
\dot{\overline A} = a\dot Z + b\overline{\dot Z},
$$
where
$a=\rho_{ZZ}+\rho_{\overline Z}\overline A_Z$ and
$b=\rho_{Z\overline Z}$.
Recall that
\begin{eqnarray*}
\dot Z =\left(
  \begin{array}{cl}
  \dot z\\
  \dot w
  \end{array}
\right)
\sim \Lambda^{-1} \left(
  \begin{array}{cl}0\\
  \varphi
  \end{array}
\right)=
\rho_z^{-1}\left(
\begin{array}{cl}
-\rho_w\\
\rho_z
\end{array}
\right) \varphi.
\end{eqnarray*}
By the Cauchy-Green formula,
$T\overline\partial\varphi=0$ and
$T\overline\partial\overline\varphi=\overline\varphi$.
Then by Lemma \ref{Lemma6},
$$
T\overline\partial\dot\lambda\vert_{b\D}
\sim
aT\overline\partial\dot Z
+bT\overline\partial\overline{\dot Z}
\sim
b\overline{\dot Z}
\sim
(\overline{\dot z},\overline{\dot w})
\left(
\begin{array}{cll}
\rho_{z\overline z}& &\rho_{w \overline z}\\
\rho_{z\overline w} & &\rho_{w \overline w}
\end{array}
\right).
$$
Then (\ref{Tdbar}) turns into
$$
\det L(\rho) \varphi \sim 0,
$$
where
\begin{eqnarray*}
\det L(\rho) = (-\rho_{\overline w}, \rho_{\overline z})
\left(
\begin{array}{cll}
\rho_{z\overline z}& &\rho_{w \overline z}\\
\rho_{z\overline w} & &\rho_{w \overline w}
\end{array}
\right)
\left(
\begin{array}{cl}
-\rho_w\\
\rho_z
\end{array}
\right)
\end{eqnarray*}
is the Levi determinant of $\rho$ with respect to $J_{st}$.
Thus $\det L(\rho)$ becomes $C^{2+\alpha}_1(b\D)$
after multiplication by any function
$\varphi \in C^{1+\alpha}(\overline\D) \cap {\cal O}(\D)$.
Hence $\det L(\rho)$ vanishes identically on the boundary
of the disc $Z_0$.

The conditions $A \circ Z_0 = 0$ and
$A_Z \circ Z_0 = 0$ imply that the Levi form of $E$
with respect to the structure $J$ coincides with the Levi
form of $E$ with respect to the structure
$J_{st}$ at every point of the boundary of the disc $Z_0$
(Lemma \ref{Leviform}).
Thus, our proposition implies that the Levi form of $E$
with respect to the structure $J$ vanishes
on the boundary of the disc $Z_0$, as desired.

This proves Proposition \ref{MainProp}
and Theorem \ref{keytheo}.

\subsection{The case of degenerate rank}
Consider the case in which the boundaries of the
pseudoholomorphic discs attached to $E$ through
a fixed point do not cover an open set in $E$.

\begin{e-pro}
\label{MainProp2}
Suppose that the boundaries of $J$-holomorphic discs
$\zeta\mapsto Z(\zeta)$ with
$Z(1) =  \left(
\begin{array}{cl}
0\\
1
\end{array}
\right)$
attached to $E$ and close to the disc
$Z_0(\zeta) =  \left(
\begin{array}{cl}
0\\
\zeta
\end{array}
\right)$, $\zeta\in\D$,
do not cover an open set in $E$.
Then for every $\zeta_0 \in b\D$, $\zeta_0 \neq 1$ there
exists a $J$-holomorphic disc attached to $E$
at the point $ \left(
\begin{array}{cl}
0\\
\zeta_0
\end{array}
\right)
$
completely contained in $E$.
\end{e-pro}

Fix a point $\zeta_0 \neq 1$ in $b\D$. For every Bishop
disc $Z(\zeta)$ close enough to $Z_0(\zeta)$
consider the evaluation map
\begin{eqnarray*}
{\cal F}_{\zeta_0}: Z \mapsto Z(\zeta_0).
\end{eqnarray*}
Then the tangent map ${\cal F}_{\zeta_0}$ of ${\cal F}$
at $Z_0$ is the map
\begin{eqnarray*}
{\cal F}'_{\zeta_0}: \dot Z \mapsto \dot Z(\zeta_0),
\end{eqnarray*}
where $\dot Z$ is an infinitesimal perturbation of $Z_0$.

Since the boundaries of discs do not cover an open subset of $E$,
we have $\rank {\cal F}'_{\zeta_0} \leq 2$ for all $\zeta_0$.
The following statement implies that
$\rank {\cal F}'_{\zeta_0} \geq 2$.

\begin{e-lemme}
\label{lemma8}
For every $w_0 \in \cc$ there exists
$\dot Z =  \left(
\begin{array}{cl}
\dot z\\
\dot w
\end{array}
\right)$
with $\dot w(\zeta_0) = w_0$.
\end{e-lemme}
\proof We recall that $A = 0$ and $A_z = 0$ on the discs $Z_0$.
Then as above we have
\begin{eqnarray*}
& &\dot Z = \Lambda^{-1}V,\\
& &V = \left(
\begin{array}{cl}
0\\
\varphi
\end{array}
\right) + R_1 \left(
\begin{array}{cl}
0\\
\varphi
\end{array}
\right) + R_2 \left(
\begin{array}{cl}
0\\
\overline\varphi
\end{array}
\right),
\end{eqnarray*}
where $\varphi$ is an arbitrary holomorphic function with
$\varphi(1) = 0$. Then
\begin{eqnarray*}
\dot w = \varphi + R_1^{22} \varphi + R_2^{22}\overline\varphi,
\end{eqnarray*}
where $R_1^{22}$, $R^{22}_2$ are $(22)$-matrix elements of
$R_1$ and $R_2$. Plugging $\zeta = \zeta_0$,
we get
\begin{eqnarray*}
\dot w(\zeta_0) = \varphi(\zeta_0)+
\int\int_{\D}a_1(\zeta)\varphi(\zeta)d\zeta\wedge d\overline\zeta
+\int\int_{\D}a_2(\zeta)\overline{\varphi(\zeta)}
d\zeta \wedge d\overline\zeta,
\end{eqnarray*}
where $a_1$ and $a_2$ are integrable in $\D$.
Assume that the rank of the map
$\varphi \mapsto \dot w(\zeta_0)$ is smaller than or equal to 1.
Then there is $c \in \cc \backslash \{ 0 \}$ such that for
every $\varphi$ we have
$Re(c \dot w(\zeta_0)) = 0$. Then for some $b_1, b_2 \in L^1(\D)$
we have
\begin{eqnarray*}
2\Re (c\varphi(\zeta_0)) + \int\int_{\D} b_1(\zeta)\varphi(\zeta)
d\zeta \wedge d\overline\zeta
+\int\int_{\D} b_2(\zeta)\overline{\varphi(\zeta)}
d\zeta \wedge d\overline\zeta = 0.
\end{eqnarray*}
Splitting into linear and antilinear parts we obtain
\begin{eqnarray*}
c\varphi(\zeta_0) + \int\int_{\D}b_1(\zeta)
\varphi(\zeta) d\zeta \wedge d\overline\zeta = 0.
\end{eqnarray*}
This implies that $c = 0$. Indeed,
take $\varphi = \psi^n$, where $\psi$ has a peak at $\zeta_0$,
that is $\psi(0) = 0$, $\psi(\zeta_0) = 1$
and $\vert \psi(\zeta) \vert < 1$ for
$\zeta \in \overline \D \backslash \{ \zeta_0 \}$.
Then passing to
the limit as $n \longrightarrow \infty$ we obtain that $c = 0$.
This contradiction proves the lemma.
\bigskip

Since the rank of the map $\varphi \mapsto \dot Z(\zeta_0)$
is equal to 2,
then there are $k_1,k_2 \in \cc$ such that
$\dot z(\zeta_0) = k_1 \dot w(\zeta_0)
+ k_2\overline{\dot w(\zeta_0)}$ for all $\varphi$.
The equality
$\dot Z = \Lambda^{-1}V$ implies
$$
\dot z = -\rho_z^{-1}\rho_w\varphi + P_1\varphi
+ P_2\overline\varphi,
$$
where $P_1$ and $P_2$ are integral operators.
Expressing $\dot z(\zeta_0)$ and $\dot w(\zeta_0)$
in terms of $\varphi$, we get
\begin{eqnarray*}
-\rho_z^{-1}\rho_w\varphi(\zeta_0) = k_1\varphi(\zeta_0)
+ k_2\overline{\varphi(\zeta_0)} +
\int\int_{\D}b_1(\zeta)\varphi(\zeta)
d\zeta \wedge d\overline\zeta
+\int\int_{\D}b_2(\zeta)\overline{\varphi(\zeta)}
d\zeta \wedge d\overline\zeta
\end{eqnarray*}
for some $b_1,b_2 \in L^1(\D)$. As in lemma \ref{lemma8} we obtain
$k_1 = -\rho_z^{-1}\rho_w\vert_{\zeta = \zeta_0}$, $k_2 = 0$.
Since $\zeta_0 \in b\D$ is arbitrary, we
have $\rho_z\dot z + \rho_w \dot w = 0$ on $b\D$ that is
\begin{eqnarray*}
\dot Z(\zeta) \in H^J_{Z(\zeta)}E, \;\;\;
\vert \zeta \vert = 1.
\end{eqnarray*}
By the hypothesis of proposition this is true for every disc
close to $Z_0$.
By the rank theorem, the image of the evaluation map
${\cal F}_{\zeta_0}$ is a
$J$-holomorphic curve in $E$. This completes the proof of
the proposition.

\subsection{Proof of Theorem \ref{MainTheorem}}

If $E$ admits a transversal Bishop  disc attached at $p$,
then the statement follows by Proposition \ref{filling}.
Suppose that there are no transversal Bishop discs.
Then by Proposition \ref{MainProp} for every disc $Z$ attached
to $E$ at $p$, the Levi form of $E$ with respect to $J$ vanishes
on $Z(b\D)$. If these discs fill an open subset $\Omega$ of $E$,
then the Levi form of $E$ vanishes on $\Omega$ identically.
Then it follows that $\Omega$ is foliated by $J$-holomorphic discs,
which holds for Levi-flat hypersurfaces in any dimension \cite{KS}.
In the case of complex dimension 2 the existence of the foliation
follows immediately from the representation (\ref{Leviformula}) for
the Levi form and the Frobenius theorem. Finally, if the boundaries
of Bishop discs do not cover an open piece of $E$, then there
exist $J$-holomorphic discs in $E$ by Proposition \ref{MainProp2}.
Thus, $E$ necessarily admits a transversal Bishop disc
which proves the theorem.

\end{document}